\RequirePackage{fix-cm}
\documentclass[amsmath,amssymb]{natureprintstyle}
\usepackage{color}
\usepackage{amsmath,amssymb,amsthm,xfrac}
\usepackage{fancyhdr}
\usepackage{psfrag}
\usepackage{empheq}
\usepackage{wrapfig}
\usepackage{array}
\usepackage{multirow}
 \usepackage{graphicx}
\usepackage{hyperref}
\usepackage{bbold}
\usepackage{enumitem}
\usepackage[table]{xcolor}
\usepackage{booktabs}

\def\dst{\displaystyle}
\def\eps{\varepsilon}

\def\p{\partial}
\def\a{\alpha}
\def\b{\beta}

\title{Flagella bending affects  macroscopic properties of bacterial  suspensions}
\author{M. Potomkin$^1$, M. Tournus$^2$, L. Berlyand $^1$, I. S. Aranson$^{3,4}$} 

\newcommand{\comMisha}[1]{\textcolor{black}{#1}}

\begin{document}
\maketitle

\begin{affiliations}
\item Pennsylvania State University, Mathematics Department, University Park, Pennsylvania 16802, USA
\item  Aix Marseille Univ, CNRS, Centrale Marseille, I2M,  Marseille, France
\item Materials Science Division, Argonne National Laboratory,
9700 South Cass Avenue, Argonne, Illinois 60439, USA
\item 
Engineering Sciences and Applied Mathematics, Northwestern University, 2145 Sheridan Road, Evanston, Illinois 60202
\end{affiliations} 


\begin{abstract} 
	To survive in harsh conditions,  motile  bacteria swim in complex environment   and respond to the surrounding flow. 
	Here  we develop a mathematical model describing how the flagella  bending \comMisha{affects} macroscopic properties of bacterial suspensions. 
	First, we show how the flagella bending contributes to the decrease of the effective viscosity observed in dilute suspension. Our results do not impose tumbling (random re-orientation) as it was done previously to explain the viscosity reduction. Second, we demonstrate a possibility of bacterium escape from the wall entrapment due to the self-induced buckling of flagella. Our results shed light on the role of flexible bacterial flagella in interactions of bacteria with shear flow and walls or obstacles. 
\end{abstract}





 \section*{Introduction }
\subsection{•}


Bacteria being among the simplest living organisms, are the most abundant  species on the planet. 
They significantly influence carbon cycling and sequestration, decomposition of biomass, transformation of contaminants in the environment. Trillions of symbiotic and pathogenic  bacteria share human body space and form microbiota.  
Behavior of bacterial suspensions is an active topic of research \cite{WuLib2000,dombrowski2004self,SokAraKesGol2007,DreDunCisGanGol2011,rusconi2014,stocker2012marine}.  
The recent discoveries include the onset of large-scale collective behavior \cite{dombrowski2004self,SokAraKesGol2007,wensink2012meso,SokAra2012}, reduction of the effective viscosity  \cite{SokAra2009,GacMinBerLinRouCle2013,LopGacDouAurCle2015}, rectification of random motion of bacteria and extraction of useful energy \cite{SokApoGrzAra2009,di2010bacterial,kaiser2014transport}, enhanced mixing in bacterial suspensions \cite{WuLib2000,SokGolFelAra2009,pushkin2014stirring,ariel2015swarming}. 


Motile bacteria utilize bundled helical flagella to propel themselves in a fluid environment. 
Bacteria use the propensity to swim to search for food (e.g. chemotaxis), colonize new territory, or escape harsh conditions. Orientation of bacteria is also affected by shear flow, leading to a variety of 
non-trivial effects, such as rheotaxis (swimming against the flow) \cite{fu2012bacterial} or depletion of bacterial concentration in shear flows \cite{rusconi2014,SokAra2016}. 
Unlike chemotaxis, i.e. drift along the concentration gradient, rheotaxis and concentration depletion are pure physical effects since no active receptor response is needed for the explanation of these phenomena. 
Elastomechanics  of the  bacteria, like bending and buckling of the flagella, could then play an important  role in the understanding of these phenomena
\cite{SonGuaSto2013}. 
A flagellum
is, typically, at least, twice longer than the bacterial
body and is flexible. Thus, flagella  bending could result in a significant effect on bacterial trajectories \cite{VogSta2012,VogSta2013,SonGuaSto2013,BruRusSonSto2015}. 
Nonlinear dynamics of rigid microswimmers  in  two-dimensional  Poiseuille flow were studied in Refs. \cite{ZotSta2012,ZotSta2013}. It was shown that the swimmers initially located away from channel walls exhibit a stable periodic motion around the centerline of the flow. Role of bacteria motility on zipping of individual flagellar filaments  and formation of bacteria flagella bundle was investigated in Ref. \cite{tapan2015}.  
However, it was poorly understood how the flagellum can affect the bacterial dynamics due to bending in response to the external shear flow or due to collision with the wall or obstacle. A model of a swimmer with flexible flagella in two fundamental shear flows, either planar shear or the Poiseuille flow in long channels, has been introduced in our previous work \cite{TouKirBerAra2015}. 
A variety of surprising  effects was discovered. For example,  depending on the bending stiffness of the flagellum, the swimmer may migrate towards the center or exhibits periodic motion. 
This paper significantly extends and advances our  results obtained in Ref.  \cite{TouKirBerAra2015}. 
Here we succeed  to tackle two new important problems associated with the bacterial dynamics in shear flows. 
We show that flexibility of the bacterial flagella (i) contributes to the reduction of the effective viscosity and (ii) assists bacteria escaping entrapment near solid walls. 
Our results provide insight how microswimmers interact with external shear flow and with obstacles, realized, for example, in microfluidic devices or {\it in vivo}.


The first part of this work is motivated  by the experimental observation in Ref. \cite{SokAra2009,GacMinBerLinRouCle2013} on the decrease of the effective viscosity of an
active suspension of {\it B. subtilis}, in particular in the dilute regime, that is the volume fraction of bacteria is less than 1\%. This result has been recently extended in Ref. \cite{LopGacDouAurCle2015} where a suspension of {\it E. coli} exhibited properties reminiscent that of a super-fluid: persistent flow and zero (or even negative)  apparent viscosity. 
This is a hallmark of active matter: chemical energy stored in nutrient is turned into mechanical energy which is then used to counter-balance the viscous dissipation. 
Suspensions of active (self-propelled) swimmers representing bacteria were studied in Refs. \cite{HaiAraBerKar2008,HaiSokAraBerKar2009,SokAra2016} and Ref. \cite{Sai2010} with the primary goal to identify a mechanism resulting in the decrease of effective viscosity in a dilute regime. The works \cite{HaiAraBerKar2008,HaiSokAraBerKar2009,Sai2010} require bacteria to tumble (randomly change direction  characterized by some tumbling rate or effective rotational diffusion $D_{\rm r}$). Nevertheless, the strain of {\it B. subtilis} used in Ref. \cite{SokAra2009}  tumbles rarely, i.e. $D_{\rm r} \ll 1$.   
Here we show that bacterial flagella bending  contributes to the reduction of the  effective viscosity     even in the absence of tumbling.      
We  derive an asymptotic expression  for the effective viscosity for a dilute suspension. We show that this expression is in agreement with both the numerical solution of the model and qualitatively consistent with the experimental data from Ref. \cite{SokAra2009}. 


The second and related part of the work focused on the bacterium behavior near surfaces (e.g. obstacles or walls). This problem naturally occurs in multiple setting relevant in biomedical context, e.g.  formation of biofilms, migration of bacteria along channels, e.g., catheter,  and industry (pipes clogging and biofouling). In many applications, bacteria swim in a confined container and their trajectory can be significantly affected by a nearby surface. Typically, bacteria are attracted by a no-slip surface (a wall) due to long-range hydrodynamic interactions \cite{BerTurBerLau2008}, then bacteria swim (mostly) parallel to the wall for a certain period of time. 
Eventually, bacteria can 
escape due to tumbling \cite{DreDunCisGanGol2011}  or can adhere to the wall. 
Study of behavior of flagellated swimmers near walls was initiated by Ref. \cite{Rot1963} where the accumulation of spermatozoa at the glass plates was documented. In the experimental works \cite{BerTur1990,diluzio2005escherichia} it was shown that {\it E. coli} is attracted by the wall and the straight trajectory becomes circular due to counter-rotation of bacterial body and the flagella. The tendency of bacteria to approach the wall and to increase the curvature of their trajectory was observed by numerical modeling in Ref. \cite{RamTulPha1993} where a bacterium was modeled as a sphere with helical flagella rotating with a constant angular velocity. To explain why the bacteria can swim near the wall adjacent to it for a long time, in Ref. \cite{FryForBerCum1995} authors hypothesized the presence of short-ranged forces of the  van der Waals type. However, in Ref. \cite{VigForWagTam2002} by combining theory and experiment it was shown that van der Waals forces cannot be responsible for parallel swimming of the bacteria near wall.
Instead authors proposed to extend the model from Ref. \cite{RamTulPha1993} for a non-spherical bacterium body, and showed that bacteria may be kept at the wall by the additional torque caused be the non-sphericity.   
In addition to hydrodynamic attraction, 
bacteria can eventually re-orient themselves  and swim away from the wall (escape). In Ref. \cite{DreDunCisGanGol2011} the time needed for bacteria to escape was estimated theoretically provided that rotational diffusion (for example, due to tumbling) is introduced. Authors \cite{DreDunCisGanGol2011} noted that even if bacteria do not tumble and are too large to be affected by thermal effects, the rotational diffusion can be assigned with a significant value due to noise in the swimming mechanism, whose essential constituent is flagellum dynamics.
Here we consider 
how  a bacterium that being initially entrapped and immobilized at a wall can escape  
 exploiting its flagella flexibility. Such an entrapment  may also naturally 
happen when the suspending liquid is anisotropic, e.g. lyotropic liquid crystal 
\cite{MusTriRoyArnWeiAbb2015,zhou2014living}. In this situation bacteria are swimming predominantly parallel to the average molecular orientation, i.e. liquid crystal director. 
In the case when the liquid crystal director is anchored perpendicular to the confining wall (homeotropic alignment),
 bacteria are forced to be aligned perpendicular to the wall and become trapped \cite{MusTriRoyArnWeiAbb2015,ZhoSokAraLav2016}.
 When the motility of bacteria is increased (by adding the oxygen),  the bacterium may turn parallel to the wall due to the torque coming from the wall and the fact that forces which kept  bacteria immobilized are  small in compare to the self-propulsion (weak surface anchoring of the liquid crystal molecules). We show that a bacterium with rigid flagellum swims along the wall, so it stays essentially entrapped. In contrast, we show that  a bacterium with flexible flagellum may rotate by  an angle larger  than $\sfrac{\pi}{2}$ and escape. \comMisha{This ability to escape reduces effects of bacteria on macroscopic properties of the suspension locally near the wall (due to decrease of bacterial concentration).}

 \comMisha{{\bf Model.}
In this work, we use a mathematical model in which a swimmer is evolving independently of the others. The model is referred to below as MMFS (the mathematical model of flagellated swimmer).
The underlying physical assumptions of this model are the following:
{\it (i)}	the two-dimensional swimmer is composed of a rigid ellipse (body) and a flexible one-dimensional segment (flagellum); the flagellum is rigidly attached to the body (clamped);
{\it (ii)}  elastic and propulsion forces on the flagellum generate the thrust force which balances the drag force and leads to the motion of the swimmer;
{\it (iii)}	a propulsion force is uniformly distributed along the flagellum; 
{\it (iv)}	background flow is not modified by the flagellum. 
The shape of the body (which is an ellipse) is described by parameter $\beta= \ell^2/(\ell^2+d^2)=1/2(1-\epsilon^2)$ where $\ell$ and $d$ are major and minor axes, respectively, and $\epsilon$ is the eccentricity of the ellipse. Small $\beta$ corresponds to rod-like bodies, and $\beta=1/2$ corresponds to spheres. The length of the flagellum,  denoted by $L$, is assumed to be constant. The full list of parameters in the model as well as their typical values can be found in Table~\ref{table_2}.}

\comMisha{We stress here that a bacterial flagellum is a flexible helical filament which exhibits propeller-like motion by rotating around its helical axis and these rotations generate the propulsion force. In MMFS, this corresponds to that, according to (ii) above, the thrust force has two separate components: due to flagellum bending (the elastic force) and due to propulsion mechanism (the propulsion force). Moreover, a flagellum is modeled as a 1D (curved) segment in a plane with no helical structure, since the propulsion force already takes into account the helical structure of the flagellum and its axial rotation. }

\comMisha{Given the initial state of the swimmer (orientation of the body and the shape of the flagellum at time $t=0$),
 MMFS entirely determines the state of the swimmer for all times $t>0$.
The unknown quantities of MMFS  are the orientation 
of the body $\theta_0(t)$, the elastic stress of the flagellum $Q(s,t)$, and the tangential angle $\theta(s,t)$ of the flagellum at the point corresponding to arc length parameter $s$; $s=0$ is at the rigid interface body/flagellum and $s=L$ is at the free end of the flagellum. 
Using basic geometric formula, given $Q(s,t)$, $\theta(s,t)$, and $\theta_0(t)$, one can recover the trajectory of the swimmer, the shape and the location of every point of the flagellum $X(s,t) = (x(s,t),y(s,t))$.
MMFS requires solving a coupled system  of an ordinary differential equation for $\theta_0(t)$ and partial differential equations for $\theta(s,t)$ and $Q(s,t)$.
The system is presented in electronic supplementary material; details on its derivation can be found in Methods.
We analyze this system both numerically and using asymptotic expansions in the regime where the flagellum is almost rigid.
}



\section*{Results}

\subsection{Effective viscosity of a dilute suspension of flagellated swimmers.}
\subsection{\it  A general formula for effective viscosity.}
Effective viscosity can be understood as a measure of the total shear stress of a suspension induced by a prescribed shear flow. 
In the context of bacterial suspensions, 
the stress resulting from the applied strain is due to the intrinsic resistance of  suspending  fluid and due to 
the stress created by the microswimmers (bacteria).
In the dilute regime (small concentration),  interactions between bacteria are negligible.  Therefore,  the superposition principle applies, that is the contribution to the total stress from all bacteria is the sum of the individual contributions. Moreover, due to their large number, each bacterium's contribution may be approximated in the sum by its expected value (taking a continuum limit).
\comMisha{The dilute framework enables us to use  MMFS to derive macroscopic properties of the suspension.}
Then the formula for the effective viscosity $\eta_{\text{eff}}$ in a
linear planar shear background flow of strain rate $\dot \gamma$ becomes \cite{Bat1967,KimKar1991,HaiAraBerKar2008,HaiSokAraBerKar2009,RyaHaiBerZieAra2011}  
\begin{equation}\label{basic_formula_for_ev}
\eta_{\text{eff}}=\eta_{0}+\sum_{i=1}^{n}\eta_{\text{bact},i}\approx \eta_{0}+n\int_0^{2\pi} \frac{\Sigma_{12}(\theta_0)+\Sigma_{21} (\theta_0)}{2\dot \gamma} P(\theta_0) d\theta_0,
\end{equation}
\comMisha{where $\eta_{0}$  is the viscosity of the suspending fluid, $n=\Phi V_L$ is the number of particles in the volume $V_L$ occupied by a suspension, $\Phi$ is the number density of bacteria,} and the integral in the right hand side of Eq. \eqref{basic_formula_for_ev} is the expected value of the contribution to the effective viscosity $\eta_{\text{bact},i}$ of the $i$th bacterium. The effective viscosity $\eta_{\text{bact},i}$ is the ratio between the anti-diagonal components $\Sigma_{12}$ and $\Sigma_{21}$ of the stress tensor $\Sigma$ (induced by the bacterium) 
 and the shear rate $\dot \gamma$. Here we assume that $\Sigma_{12}$ and $\Sigma_{21}$ are only determined by the angle orientation $\theta_0$ of the bacterium. Thus, in order to compute the expected value of $\Sigma_{12}$ and $\Sigma_{21}$, finding the distribution of orientation angles $P(\theta_0)$ is necessary. 
 
 \comMisha{MMFS is based on balance of forces and torques exerted by the swimmer on its rigid body surface and at each point of the flexible flagellum and the fluid drag forces and torques, respectively. Moreover, in MMFS the sum of forces exerted by the swimmer is zero. This is similar to the force-dipole model of a swimmer\cite{DreDunCisGanGol2011,RyaHaiBerZieAra2011} where the sum of the force that pushes the body in the fluid and the force that perturbs the fluid due to propulsion mechanism in the flagellum (represented in the force-dipole model by a point force exerted behind the body) is zero. The key difference between MMFS and the force-dipole model is that the sum of all torques exerted by the swimmer in MMFS is not necessarily zero, whereas for force-dipoles this sum is trivially zero. In particular, the force-dipole can not rotate if no external torque is exerted (a non-zero background flow, interactions with other swimmers, external magnetic field, etc.), while the flagellated swimmer may rotate if the flagellum is bent. The fact that the fluid balances a non-zero total torque exerted by the swimmer results in that, in general, the effective stress is non-symmetric in MMFS, i.e., $\Sigma_{12}\neq \Sigma_{12}$\cite{Bat1970}. Presence of a non-zero anti-symmetric part of the effective stress due to active contribution is the special feature of active chiral fluids\cite{FurStrGriJul2012}.}


In order to find $\Sigma_{kl}$ one needs to solve the Stokes equation in the low Reynolds number regime: $-\nabla_{\bf x}\cdot\Sigma({\bf x}) = F_{\text{bact}}({\bf x})$, where $\Sigma$ is the fluid stress tensor, and $F_{\text{bact}}$ is the bulk force due to the presence of the bacterium (\comMisha{the thrust force}). Solving this equation is impractical due to the large domain of integration compare to the bacterium size. In a simpler model of a bacterium with rigid flagellum \cite{RyaHaiBerZieAra2011}, each bacterium was approximated by a force dipole \cite{DreDunCisGanGol2011} and the explicit expression for $\Sigma$ is well-known in this case. Here our goal is to capture the effects coming from bending of elastic flagella in shear flow, so an approximation by the force dipole would oversimplify the consideration and would lead to zero net contribution to the effective viscosity. Instead, we use the Kirkwood approximation  for the stress tensor \cite{KirAue1951,DoiEdw1988,ZieAra2008}   
\begin{equation}
\label{kirkwood_formula_general}
\Sigma_{kl}=\frac{1}{V_L  }\int (F_{\text{bact}}({\bf x}))_k({\bf x}-{\bf x}_{c})_l d{\bf x},
\end{equation} 
where $F_{\text{bact}}({\bf x})$ is non-zero in a small neighborhood of the center of mass of the bacterium ${\bf x}_{c}$, and $V_L$ is the volume occupied by the fluid. 
The Kirkwood approximation  can be also interpreted as the second term in the multi-pole expansion \cite{KimKar1991}.  

\comMisha{In the context of MMFS (defined in the end of Introduction),}
$F_{\text{bact}}({\bf x})$ is the sum of two forces distributed over the flagellum: $(i)$ the uniform propulsion force $F_p \boldsymbol{\tau}(s)$ directed along the unit tangent vector $\boldsymbol{\tau}(s)$ with the magnitude $F_p$ and $(ii)$ the elastic force $Q(s)=\Lambda (s){\boldsymbol{\tau}}(s) + N (s){\boldsymbol{n}}(s)$ ($\Lambda$ and $N$ are tangent and normal components of $Q$, respectively). 
\comMisha{In the following, it will be convenient to separate contributions coming from propulsion and from elasticity 
for the components of the stress tensor: $\Sigma_{kl}=\Sigma_{kl}^{\text{propulsion}}+\Sigma_{kl}^{\text{elastic}}$, where according to the formula \eqref{kirkwood_formula_general}
\begin{equation}\label{sigmas}
\begin{aligned}
 \Sigma_{kl}^{\text{propulsion}} &= \frac{1}{V_L}\int_0^L F_p \boldsymbol{\tau}_k(s) (X_l(s) -X_l(0)) ds , \\
 \Sigma_{kl}^{\text{elastic}} &= \frac{1}{V_L }\int_0^L  \frac{\p Q_k}{\p s}(s)(X_l(s) -X_l(0))ds,
\end{aligned}
\end{equation}
and analogously for the effective viscosity $\eta_{\text{eff}}$:
\begin{equation}\label{formula_for_eta_contributions}
\dfrac{\eta_{\text{eff}}-\eta_{0}}{\eta_{0}}=\eta_{\text{propulsion}}+\eta_{\text{elastic}},
\end{equation}
where each of the two terms is computed via \eqref{sigmas} and \eqref{basic_formula_for_ev}. 
}


\comMisha{Terms in the right hand side of \eqref{formula_for_eta_contributions} take into account the effect of the flagellum. This means that the contribution to the effective viscosity due to presence of rigid ellipsoidal bodies in the fluid is included in $\eta_0$. In other words, $\eta_0$ is the effective viscosity of the dilute suspension of rigid ellipsoids and $\eta_0\approx\eta_{\text{fluid}}\left[1+\nu \Phi\right]$, where $\eta_{\text{fluid}}$ is the viscosity of water, $\Phi$ is the number density. The formula for the coefficient $\nu$ is well-known: for spheres $\nu=2.5$ (the Einstein's formula \cite{Ein1906}), for ellipsoids the formula for $\nu$ was obtained by Jeffery (formulas (62) and (64) in Ref.~\cite{Jef1922}).}

  
  \medskip
\subsection{ \it Asymptotic results for large bending stiffness of flagellum.}
We present here  our results on computations of $\eta_{\text{propulsion}}$ and  $\eta_{\text{elastic}}$ as functions of the following geometrical and physical dimensionless parameters: shape parameter $\beta$ (describes shape of the bacterium body; $\beta = 0$ for needles and $\sfrac12$ for spheres), ratio of bacterial body length to the flagella length $r=\ell/L$ ($\ell$ and $L$ are the body and the flagellum length, respectively), and the compound dimensionless parameter characterizing ratio of drag  force to elastic force $\eps=L^4 \dot \gamma \comMisha{\zeta_b}/K_b$ ($\comMisha{\zeta_b}$ is the drag coefficient and $K_b$ is the bending stiffness of the flagellum). 
  We use two scale asymptotic expansions in small $\eps$ (stiff flagella) to establish explicit expression  for the effective viscosity $\eta_{\text{eff}}$. 
Note that for fixed values of $\comMisha{\zeta_b}$, $L$ and $\dot \gamma$, taking $\eps$ small is equivalent to the limit as flagellum is nearly rigid.  That implies that  the bending stiffness of the flagellum $K_b$ is large (the reader should not be confused by the fact that the typical $K_b$ we use for bacteria and call it ``large" is of the order $10^{-23}$ N$\cdot\text{m}^2$; after nondimensionalization $K_b$ is replaced by $\eps^{-1}$, for details see \comMisha{electronic supplementary material, section 2}). 
Necessity of two scales in the asymptotic expansions in the "rigid" limit is explained by the two different time scales for dynamics of the smoothly translating bacterial body and rapidly oscillating flagellum.

\comMisha{The two scale asymptotic expansion methods for equations of MMFS (see supplementary material, section 2, for a description of the method)} is used to derive the following asymptotic expression for the tangential and normal components of the elastic stress: 
\comMisha{ $\Lambda(s)=p_{\Lambda}(s)\sin 2\theta_0\!-\!\frac{F_p(s-L)}{1+k_r},$ $N(s) = p_{N}(s)\cos 2\theta_0$.  }
\comMisha{Polynomials $p_{\Lambda}(s)$ and $p_{N}(s)$ are of the second order with respect to arc length $s$ with coefficients proportional to $\comMisha{\zeta_b}$ and $\dot{\gamma}$ and they also depend on shape parameter $\beta$, flagellum length $L$, body length $\ell$ and drag coefficient $k_r$  (see Table \ref{table_2} with the list of parameters). The second term in the expression for the tangential component of the elastic stress is due to the propulsion force which acts in the tangential direction $\boldsymbol{\tau}$ with the strength $F_p$.     
We also found the asymptotic expression for the flagellum shape described by the slope angle: $\theta(s) = \theta_0 + \eps \;p_{\theta}(s)\,\cos 2 \theta_0,$ where $p_{\theta}(s)$ is a polynomial of the fourth order with respect to arc length $s $ and coefficients depending on $\beta$, $\ell$, $L$, $k_r$. 
}
 \comMisha{Details of derivation of expressions for $\Lambda$, $N$ and $\theta$ with explicit formulas for coefficients of polynomials $p_\Lambda$, $p_N$, and $p_{\theta}$ 
can be found in the electronic supplementary material.}   

{The distribution of orientation angles  $P(\theta_0)$ from Eq. \eqref{basic_formula_for_ev} is in general a function of both angle of the body $\theta_0$ and time $t$,  and satisfies the Liouville continuity equation  
	\begin{equation}\label{Liouville}
	\left\{
	\begin{aligned}
	& \dfrac{\p }{\p t}P(\theta_0,t)  + \dfrac{\p }{\p \theta_0}\Big[ \frac{1}{\zeta_r}\left(\text{T}_{\text{shear}}+\text{T}_{\text{flagellum}}\right) P(\theta_0,t)\Big] = 0, \\
	& \dst\int_0^{2 \pi}  P(\theta_0,t)  \; d \theta_0 = 1,
	\end{aligned}
	\right.
	\end{equation}
	where $\text{T}_{\text{shear}}$ and $\text{T}_{\text{flagellum}}$ are torques exerted on the ellipsoidal body of the bacterium by the background shear flow and by the flagellum, respectively. \comMisha{Parameter $\zeta_r$ is the rotation drag coefficient, and  $\frac{1}{\zeta_r}\left(\text{T}_{\text{shear}}+\text{T}_{\text{flagellum}}\right)$ is the angular velocity of the body caused by shear and flagellum}. It is well-known\cite{KimKar1991} that $\text{T}_{\text{shear}}$  can be explicitly written as a function of $\theta_0$:
	 \begin{eqnarray}\nonumber
	\text{T}_{\text{shear}}&=&\comMisha{-\dot\gamma \zeta_r((1-\beta)\sin^2\theta_0+\beta\cos^2\theta_0)}\\&=&\comMisha{-\frac{\dot\gamma \zeta_r}{2}\left(1-(1-2\beta) \cos 2\theta_0\right).}\label{def_of_T_shear}
	 \end{eqnarray} 
	 Equality $\zeta_r\frac{d\theta_0}{dt}=\text{T}_{\text{shear}}(\theta_0)$ is known as the Jeffery equation for rotating \comMisha{ellipses} in the shear flow \cite{Jef1922,KimKar1991}. In order to compute
	\comMisha{$ 
	 \text{T}_{\text{flagellum}} =  \dfrac{\ell}{2}\;N|_{s=0},
	$}
	 one needs to solve the elasticity equations for the flagellum.
	 \comMisha{
	 	However, using the asymptotic method it is possible for $\eps\ll 1$ to represent $N|_{s=0}$ (and, thus, $\text{T}_{\text{flagellum}}$) as a function of $\theta_0$, which turns Eq. \eqref{Liouville} in a closed form.} \comMisha{The resulting equation is the same as the Jeffery equation for ellipses with the effective shape parameter $b=\dfrac{r\beta}{1+2r}$ in place of $\beta$. In other words, an ellipse with the rigid flagellum has same trajectories as the more prolate ellipse with no flagellum.} The equilibrium distribution which satisfies Eq. \eqref{Liouville} for $\eps \ll 1$ is given by 
	\begin{equation}
	\label{distrib_p_theta}
	P(\theta_0)= \frac{q}{2\pi}\;\dfrac{1}{1-(1-2b)\cos (2\theta_0)}, 
	\end{equation} 
where $q=\sqrt{1-(1-2b)^2}$; \comMisha{constant $q/2\pi$ is introduced, so $P(\theta_0)$ satisfies the normalization condition in \eqref{Liouville}}.
	Since the effective viscosity should be a property of the suspension independent of time and all solutions of Eq. \eqref{Liouville} for $\eps\ll 1$ converges to the equilibrium distribution (if one assumes a small rotational diffusion), we use $P(\theta_0)$ from Eq.  \eqref{distrib_p_theta} when apply Eq. \eqref{basic_formula_for_ev}.

%
\comMisha{Substituting the asymptotic expansions into formula \eqref{basic_formula_for_ev},} the effective viscosity of the dilute suspension of flagellated swimmers is expressed as
\begin{equation}\label{actual_formula_for_ev}
\dfrac{\eta_{\text{eff}}-\eta_0}{\eta_0} =
	\Phi\frac{ L^3}{\eta_0}Z_{\text{elastic}}(\beta,r)
	-\Phi \eps \dfrac{ F_p  L^2}{ \eta_0\dot \gamma}Z_{\text{prop}}(\beta,r), 
\end{equation}
where $\Phi$ is the number density of bacteria in the suspension and expressions for elastic and propulsion contribution coefficients $Z_{\text{elastic}}$ and $ Z_{\text{prop}}$ can be found in the electronic supplementary material, {section 2.5}; \comMisha{both $Z_{\text{elastic}}$ and $ Z_{\text{prop}}$ are positive}. We point out here that Eq. \eqref{actual_formula_for_ev} implies  that the change of the effective viscosity is obtained by the  interplay between elastic and propulsion contributions. 
Namely, whatever the parameters $\beta$ and $r$ are, the propulsion  decreases the viscosity, whereas the elastic part of the stress tends to increase 
the viscosity (see Figure \ref{fig:num_asympt}).
For small $r$ (i.e. \comMisha{long} flagella), the propulsion contribution coefficient $Z_{\text{prop}}$ behaves as $ 1/r^5$ whereas the elastic contribution coefficient $Z_{\text{elastic}}$ behaves as $1/(r^{2}\log(r))$.
This implies that for $r$ small enough, the propulsion should dominate elasticity.

\begin{figure*}[t]
	\centering
	\begin{tabular}{cc}
		\includegraphics[width=0.45\textwidth]{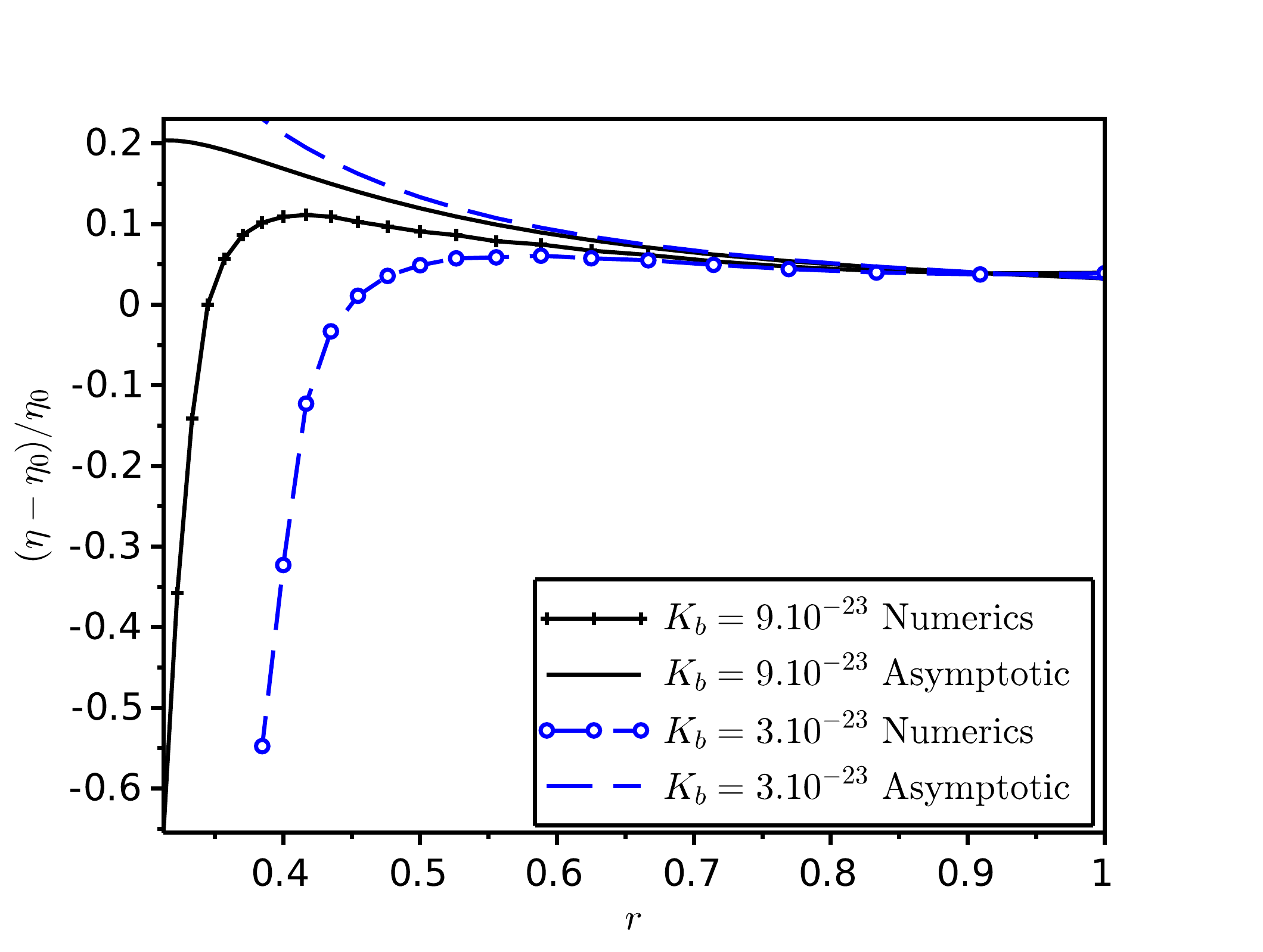}&
\includegraphics[width=0.45\textwidth]{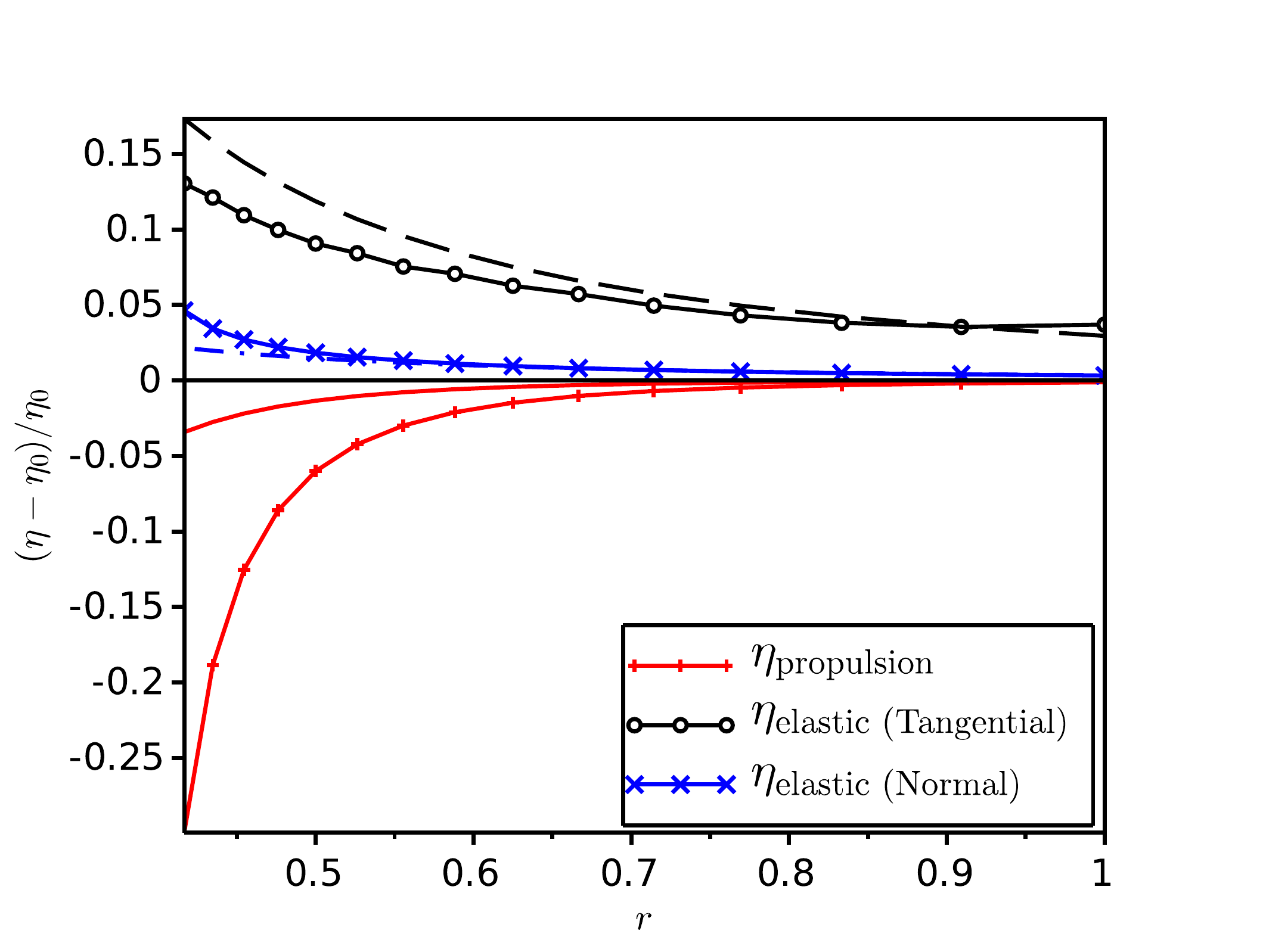}\\ (a)&(b)\\
	\includegraphics[width=.48 \textwidth]{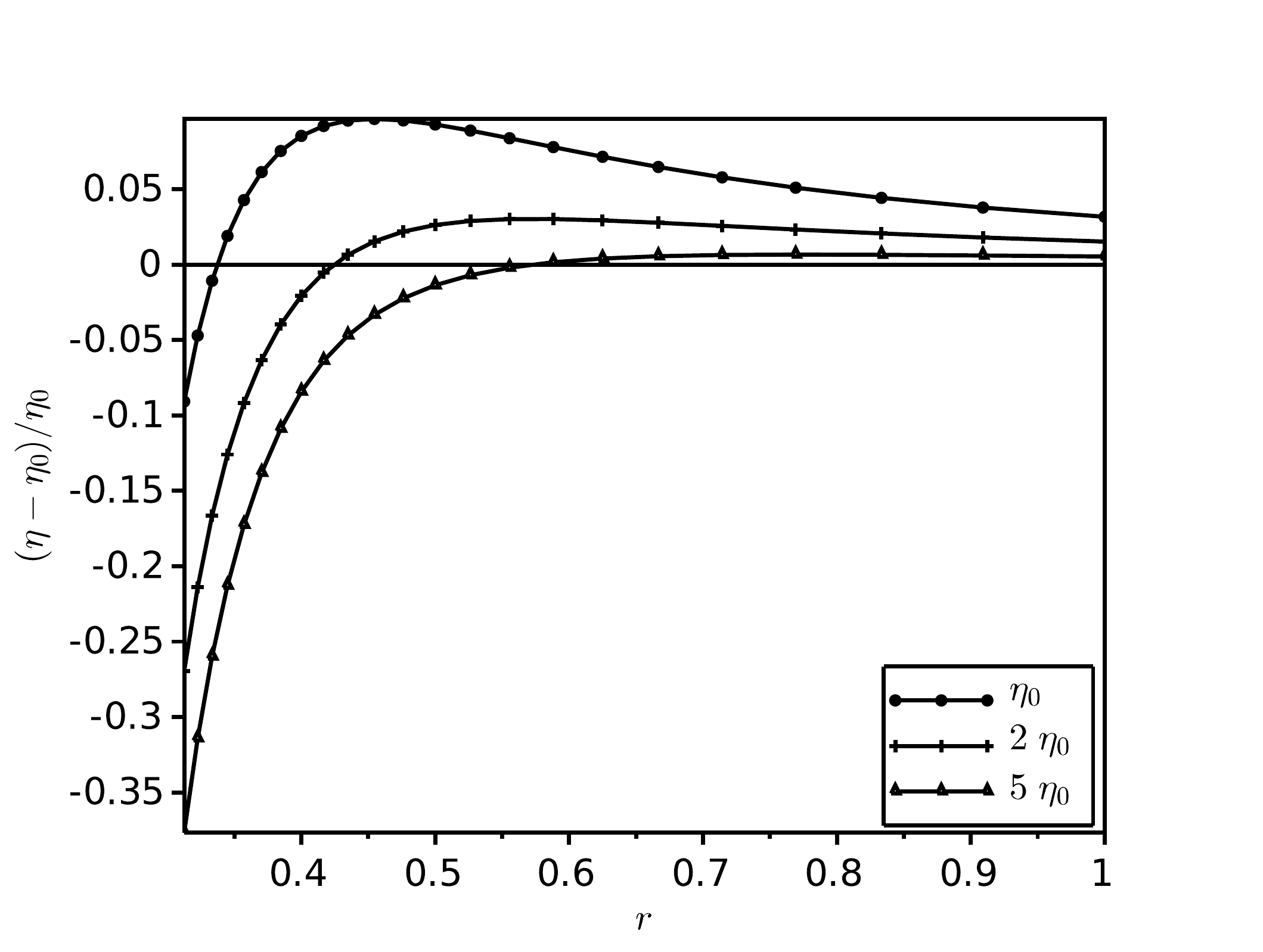}&
	\includegraphics[width=.48 \textwidth]{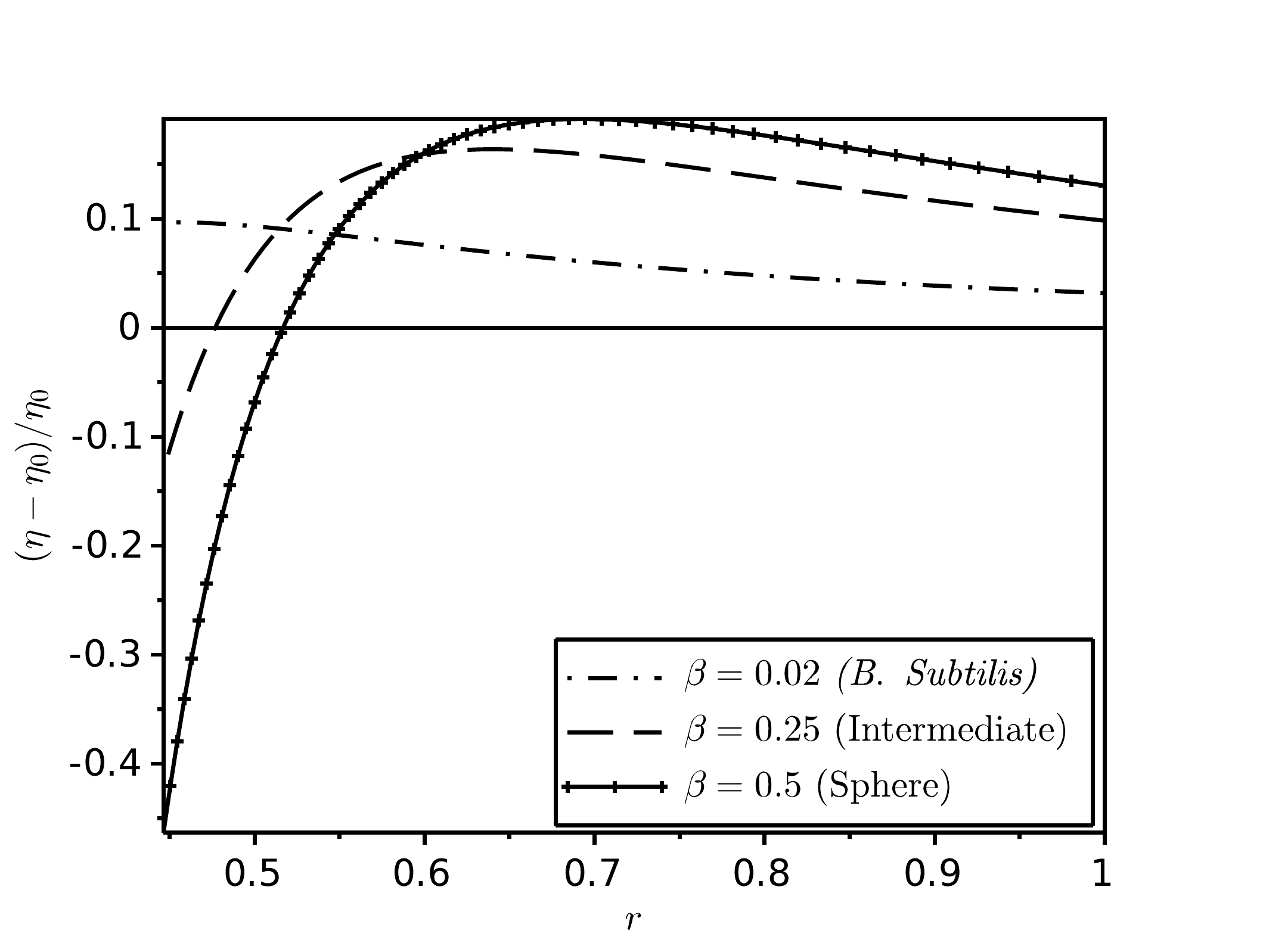}\\ (c)&(d)
	\end{tabular}
	\caption{{\bf Effective viscosity as a function of model parameters}. (a) Comparison of effective viscosity change obtained by  two scale asymptotic expansions and numerical simulations  for $K_b=3\cdot 10^{-23} \;\text{N}\cdot\text{m}^2$ (black) and $K_b$$=$$9\cdot 10^{-23}$ $\text{N}\cdot\text{m}^2$ (blue) and various $r$; (b) Contributions from propulsion and elasticity (both tangential and normal components) for $K_b=3\cdot 10^{-23} \;\text{N}\cdot\text{m}^2$ are shown.
		 (c) Effective viscosity $\eta_{\text{eff}}$ for various fluid viscosities $\eta_0$.  (d) Effective viscosity vs the body shape constants $\beta$.
		}  
	\label{fig:num_asympt}
\end{figure*}

In order to present quantitative results on the effective viscosity obtained by two scale asymptotic expansions, consider a dilute suspension with the volume fraction of bacteria about 1\%. 
The total force generated by the flagellum is about 10 pN, and the 
typical size of the flagellum is about 10 $\mu$m, which makes the propulsion strength equal to $10^{-11} \cdot 10^{-5} = 10^{-16}$ $\text{N}\cdot \text{m}$. 
The suspending  fluid is taken to be water ($\eta_0=10^{-3}$ $\text{Pa}\cdot\text{s}$). 
A realistic value of the bending stiffness of the flagellum is $K_b=3\cdot 10^{-23}$ $\text{N}\cdot\text{m}^2$ (\cite{DarBer2007}) and shear rate $\dot \gamma= 0.1\; \text{s}^{-1}$,  the parameter $\eps = 0.03$.
The values of $\eta_\text{eff}$ computed by \eqref{actual_formula_for_ev} is depicted in Figure~\ref{fig:num_asympt}. 


Certain flagellated bacteria have the ability to swim through environments of relatively high viscosity \cite{GreCan1977,SchDoe1974,Kel1974}.
Then the bacteria can maintain an almost  constant speed whatever the fluid resistance they encounter. 
In the first approximation, the velocity of the swimmer $\sim F_p /\eta_0$, then bacteria can increase their propulsion force while surrounded by a more 
viscous fluid.
In such a fluid ($\eta_0 = 5.10^{-3}$ $\text{Pa}\cdot\text{s}$ and $F_p = 1.5 \; \mu \text{N}\cdot \text{m}^{-1}$), we predict a decrease of viscosity for $r<0.55$ ($r=0.5$ for {\it B. subtilis}) (Figure \ref{fig:num_asympt}, c).
For higher values of the shape constant $\beta$, the decrease of effective viscosity  also occurs for shorter flagella (Figure \ref{fig:num_asympt}, (d)). 


  \medskip
\subsection{\it Numerical simulations.}
We performed computational analysis of  
\comMisha{MMFS} and  computed the effective viscosity  $\eta_{\text{eff}}$ as well as propulsion and elastic contributions $\eta_{\text{propulsion}}$ and $\eta_{\text{elastic}}$. The expected value integral in the right hand side of Eq. \eqref{basic_formula_for_ev} was approximated by the time average of its integrand. 
 
For large $K_b$ (small $\eps$), the results of numerical solution  are in a good agreement with asymptotic expression  Eq. \eqref{actual_formula_for_ev} (see Figure \ref{fig:num_asympt} and the Table \ref{table1}).
Note that the asymptotic parameter $\eps$ is proportional to $L^4$, so the agreement between numerical  and asymptotic solution is lost for small $r$ (long flagellum).
The set of values of the flagellum length $L$ for which the decrease of viscosity is observed depends on the bending stiffness $K_b$.
Table \ref{table1}  compares results of asymptotic approach and numerical solution, the threshold values of $L$, $r$ and $\eps$ needed to have a decrease of viscosity are given.
\begin{table*} 
\begin{center}
\begin{tabular}{|c|c c|c c|}
\hline
\multirow{2}{*}{} &\multicolumn{2}{c|}{ $K_b$=$3\cdot 10^{-23}\;\text{N}\cdot\text{m}^2$}  & \multicolumn{2}{c|}{ $K_b$=$9\cdot 10^{-23}\;\text{N}\cdot\text{m}^2$}\\
\cline{2-5}
 &  Asymptotics  & Numerics  & Asymptotics & Numerics  \\
 \hline 
  & $L >15 \mu \text{m}$ &  $L >11 \mu \text{m}$&  $L >22 \mu\text{m}$  & $L >15 \mu\text{m}$ \\
   $\dfrac{\eta_{\text{eff}} -\eta_0}{\eta_0} < 0$    &   $r< 0.33$  &$r< 0.45$  & $r< 0.23$ &$r< 0.34$ \\ 
        &$\eps >  0.16 $ &$\eps >  0.05 $ &  $\eps>0.24$   &$\eps >  0.05 $ \\
 \hline
& $L >16 \mu\text{m}$  &$L >12 \mu\text{m}$  & $L >23 \mu\text{m}$    & $L >16 \mu\text{m}$\\ 
 $\dfrac{\eta_{\text{eff}} -\eta_0}{\eta_0} < 10\% $ & $r< 0.31$  &  $r< 0.43$  &  $r< 0.22$  & $r< 0.32$ \\
 & $\eps >  0.21 $ & $\eps >  0.06 $ & $\eps> 0.26$   &$\eps >  0.06 $\\
 \hline
\end{tabular}
\caption{Comparison of numerical solution with the asymptotic results}
\label{table1}
\end{center}
\end{table*}

For the following model parameters:   $K_b = 3\cdot 10^{-23}\;\text{N}\cdot\text{m}^2$, $\eta_0=10^{-3}\,\text{Pa}\cdot\text{s}$, $L=12\, \mu \text{m}$ and $F_p=1.5\, \mu\text{N}\cdot \text{m}^{-1}$, 
asymptotic and numerical values of $\eta_{\text{eff}}$ are in agreement with experiment \cite{SokAra2009}, i.e. we predict a decrease of effective viscosity of $\approx 10\%$
for the number density of $\Phi=5\cdot 10^9 \;\text{cm}^{-3}$ (see the first part, $\Phi<10^{9}\;\text{cm}^{-3}$, of the curve in Figure 3 in Ref. \cite{SokAra2009}).

\subsection{Flagellated swimmers can escape from the wall.}

 
Here we consider how flagella flexibility assists the bacterium to escape from the wall. 
A swimmer can be entrapped by a wall such that its orientation is perpendicular to the wall. This kind of entrapment may happen,  for example,  in lyotropic (water soluble) nematic liquid crystal with the homeotropic surface  anchoring \cite{MusTriRoyArnWeiAbb2015,ZhoSokAraLav2016} (the liquid crystal director is perpendicular to the wall): 
since the bacteria tend to align with the nematic director, they eventually become perpendicular to the wall (for simplicity we neglect here the effects associated with the anisotropic elastic and viscous torques exerted by the liquid crystal on a bacterium). Moreover, motile bacteria would hit the wall. 
However, due to flagella  rotation and bending, this perpendicular alignment may become unstable. 

Settings of the problem are as follows. Bacterium's body initially has the orientation $\theta_0=\pi$, that is, the body is oriented horizontally, pointing to the right, at the vertical wall $x=0$. Flagellum is initially  slightly perturbed  from a straight configuration (while unstable, perfectly straight flagellum will lead to no motion). Bacterium's body experiences  three torques: $(i)$ due to the flagellum, applied at the point of its attachment to the body, $(ii)$ due to the wall, applied at the point of touching the wall (if the body does not touch the wall, then this torque is 0), and $(iii)$ due to the surrounding viscous fluid (see Figure \ref{fig:escape}, (a)). 

\begin{figure*}[]
	\begin{center}
		\begin{tabular}{cc}
		\includegraphics[width=0.48\textwidth]{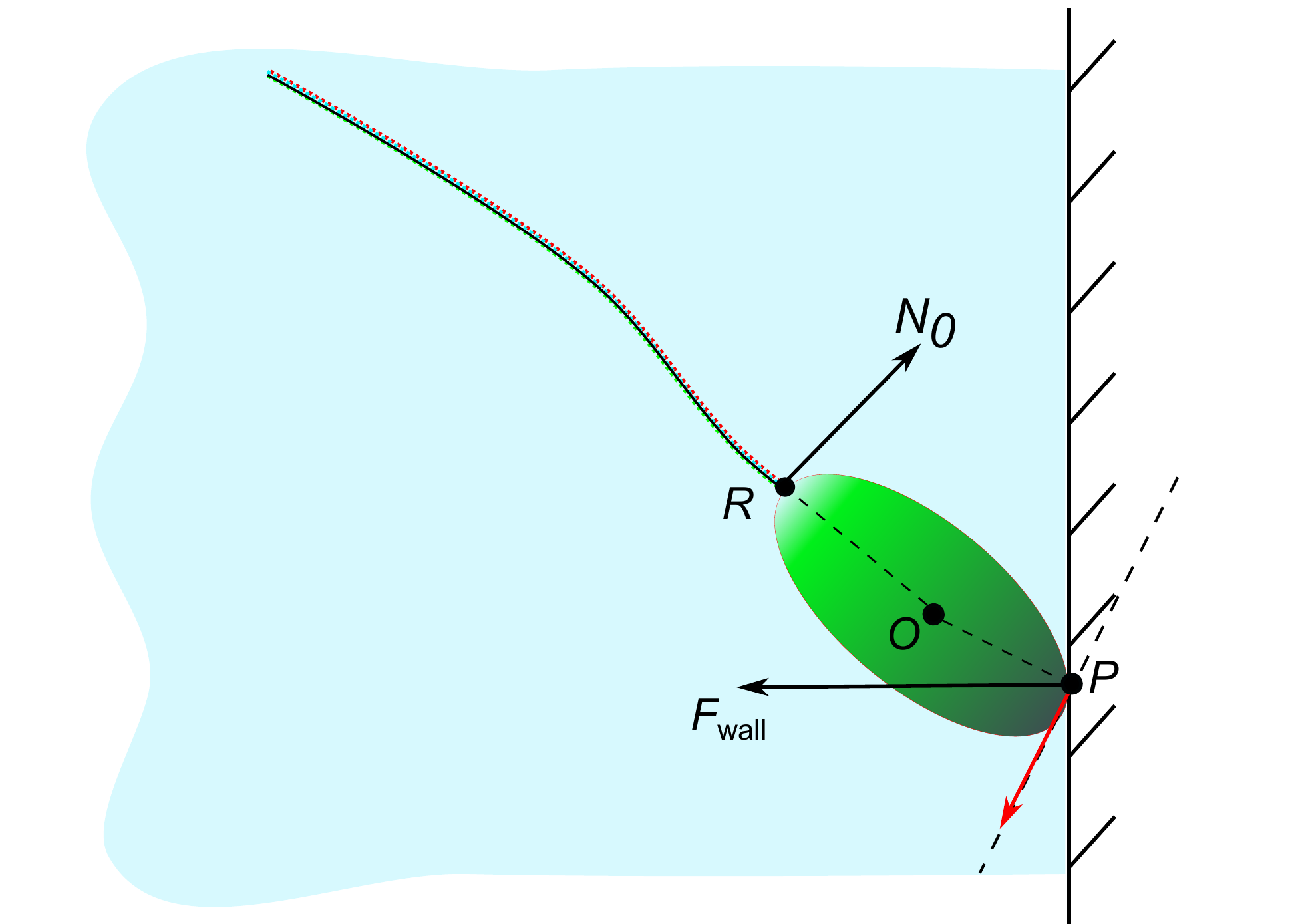}& 
		\includegraphics[width=0.46\textwidth]{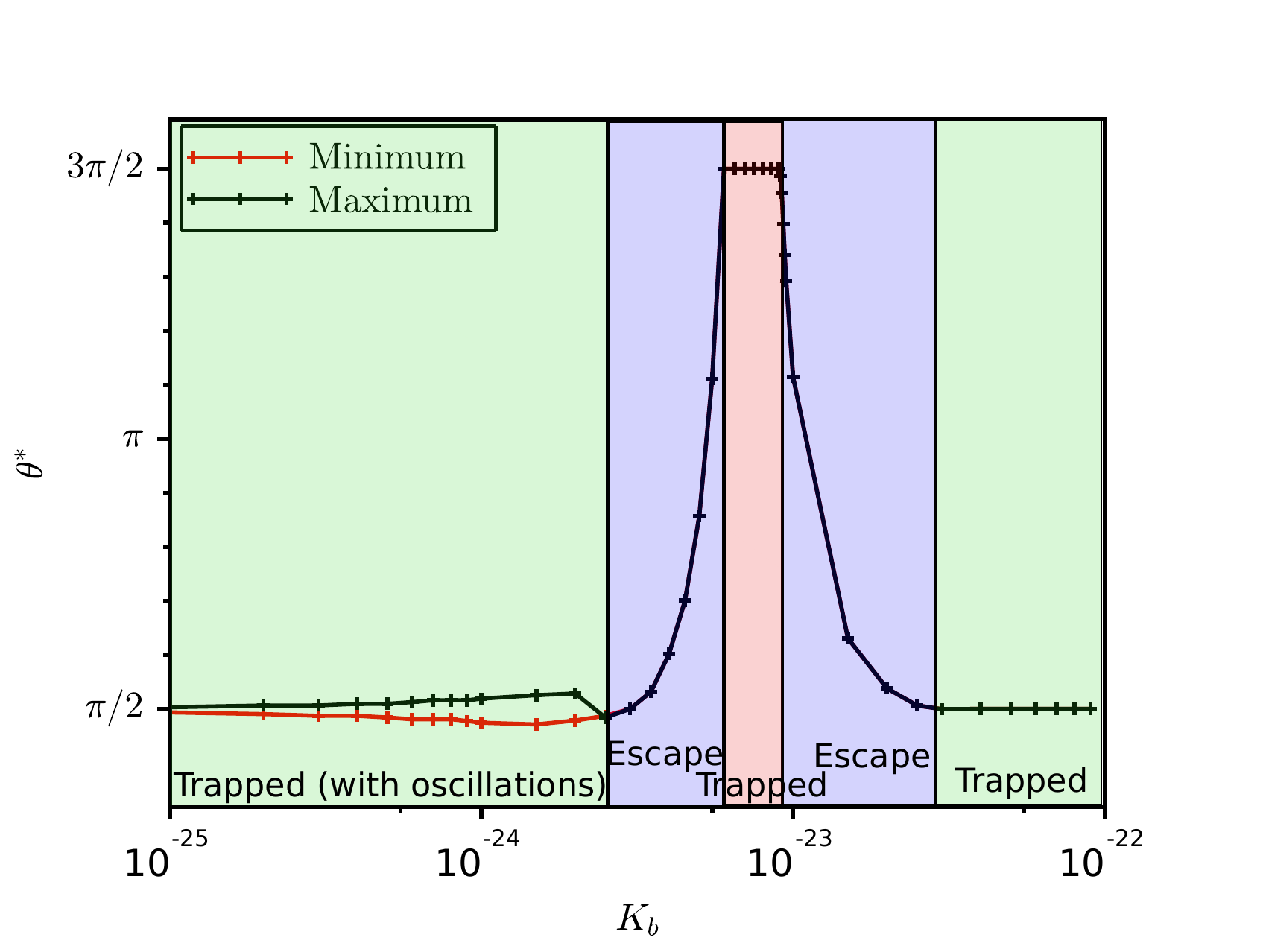}\\(a)&(b)\\  
		\includegraphics[width=0.48\textwidth]{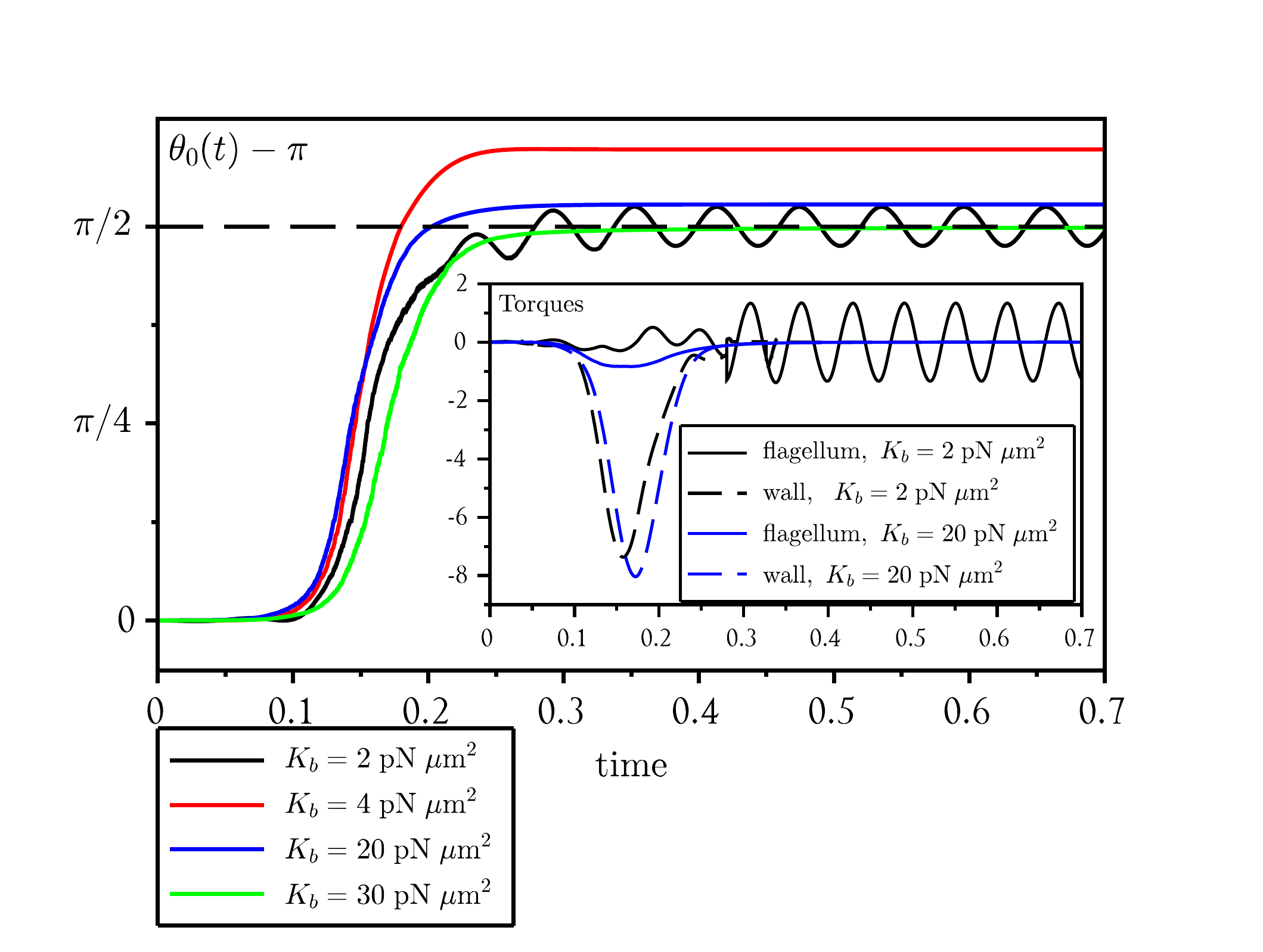}&
		\includegraphics[width=0.48\textwidth]{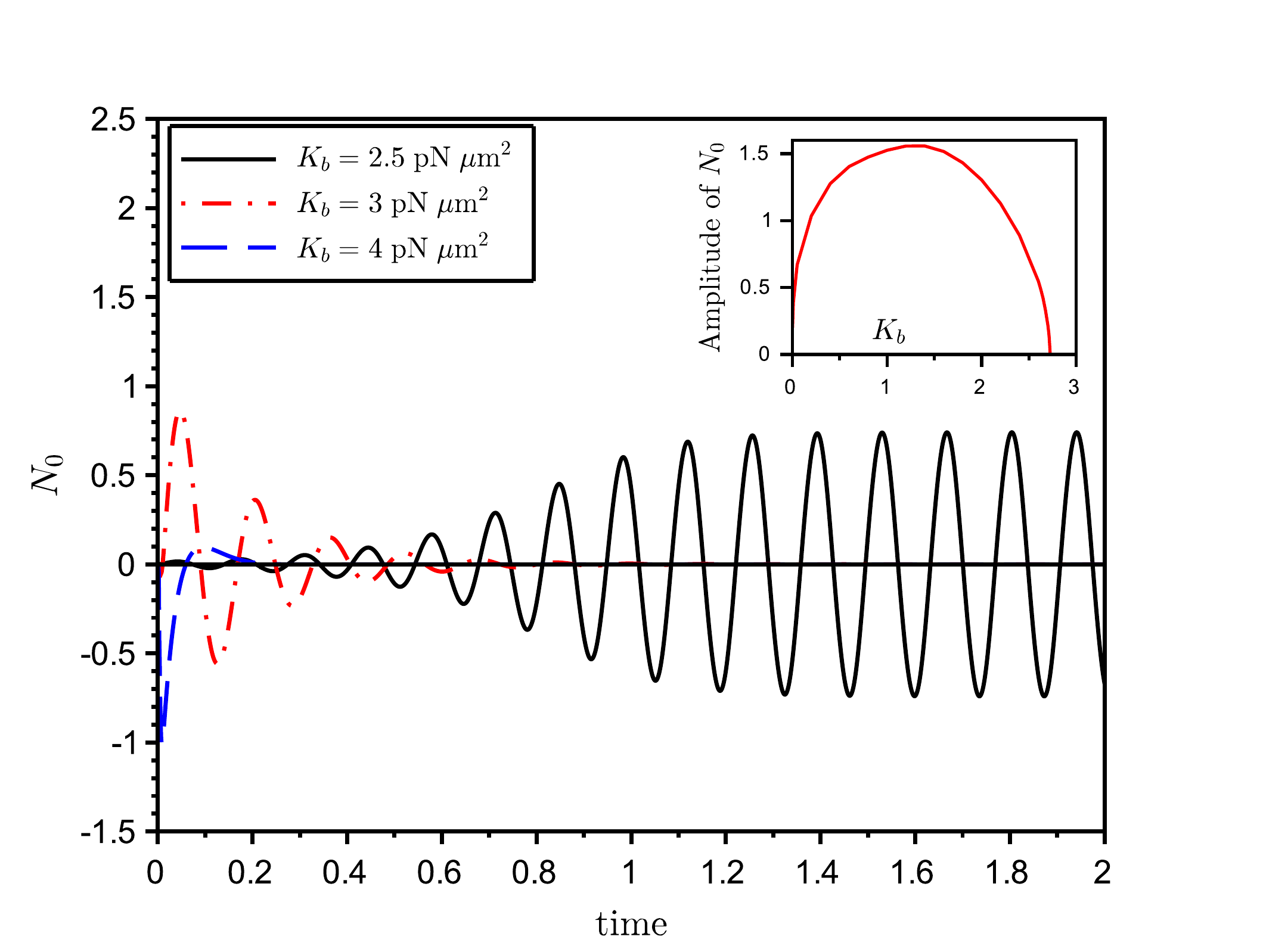}\\
		(c)&(d)
		\end{tabular}
		\caption{{\bf Escape from the wall}. 
		(a)  Bacteria at the wall (sketch). Body is an ellipsoid  centered at point $O$; the ellipsoid touches the wall at point $P$; flagellum is rigidly attached to the body at point $R$. Three forces act on the body: the normal force coming from flagellum $N_0$, the wall reaction $F_{\text{wall}}$, and the fluid over the boundary of the ellipsoid. Red arrow stands for the component of $F_{\text{wall}}$ contributing to the torque. The light blue zone represents the domain with fluid. (b)  The plot demonstrates how qualitative behavior of an initially entrapped swimmer depends on flagellum bending stiffness $K_b$: the swimmer either remains trapped (red and green zones), or eventually escapes with a limiting angle from $(\frac{\pi}{2}, \frac{3 \pi}{2})$ (blue zone); $\theta^*$ (vertical axis) denotes the orientation of the swimmer for large times, $t\gg 1$. 
		(c) Evolution of the orientation angle (main plot) and torques due to flagellum and the wall (inset). A bacterium swims towards the wall and touches it at
		time $t\approx 0.1$. When the bacterium  body touches the wall, the torque due to the wall becomes non-zero.  
			(d) The plot depicts dynamics $N_0$ for various $K_b$ of the swimmer in the fluid with no obstacles and no background flow (the plot for $K_b=4\cdot 10^{-24} \;\text{N}\cdot\text{m}$ is magnified by a factor 200 for better
			visibility); three plots demonstrate that depending on $K_b$ the swimmer eventually exhibits oscillations with constant or decaying amplitude or $N_0$ converges to $0$ with no oscillations. Inset: the plot demonstrates the dependence on $K_b$ of the amplitude of $N_0$ when the swimmer exhibits oscillations.  	  
		} 
		\label{fig:escape}
		\end{center}		
\end{figure*}

\comMisha{We use MMFS with a modification to take into account the 
 additional torque when the swimmer is touching the wall.}
A variety of nontrivial swimming regimes was numerically observed depending on the values of the bending stiffness of the flagellum, $K_b$, 
\comMisha{all other parameters being fixed, their values can be found in Table 2}. 
Numerical analysis shows that qualitative behavior of the swimmer depends on the bending stiffness of the flagellum, $K_b$. If $K_b<5\cdot 10^{-24}\;\text{N}\cdot\text{m}^2$ ("soft" flagellum) the swimmer rotates and swims parallel to the wall. Thus, in this case, though the swimmer is not immobilized at the wall, it is still entrapped by the wall and cannot escape. For $5\cdot 10^{-24} \;\text{N}\cdot\text{m}^2<K_b<2.2\cdot 10^{-23}\;\text{N}\cdot\text{m}^2$, the swimmer eventually swims away from the wall, hence showing ability to escape due to the flagellum. For the large bending stiffness, $K_b>2.2\cdot 10^{-23}\;\text{N}\cdot\text{m}^2$ ("almost rigid" flagellum), as it is expected for the rigid flagellum, swimmer rotates by $\pi/2$ and then swims parallel to the wall (see Fig. \ref{fig:escape}, (b) and (c), and the electronic supplementary material,
videos S1--S4).   
The difference between "soft" and "almost rigid" cases is that in the first one the body exhibits visible oscillations, whereas in the latter, it swims straight parallel to the wall. This non-trivial qualitative behavior of the swimmer depending on $K_b$ is also observable when no obstacle is present in the fluid. Regardless initial shape the swimmer eventually either orients itself toward one direction and swims straight or exhibits periodic oscillations (see Figure  \ref{fig:escape}, (d) for the case with no background flow; more complicated dependence on $K_b$ of large time behavior of the swimmer was observed in the Poiseuille flow, see Figure 3 (e) in Ref.\cite{TouKirBerAra2015}; see also the electronic supplementary material, video S5).
These differences in qualitative behavior may serve as a basis to isolate bacteria with bending stiffness in a given range (or equivalently, different numbers of flagella since effective $K_b$ is proportional to the number of flagella).


\section*{Discussion}
\subsection{•}

We provide a heuristic explanation why the flagellum helps decrease the viscosity. For the illustration we will use the force dipole representation of a bacterium, that is, the representation by two forces of equal magnitude and opposite directions. In work \cite{DreDunCisGanGol2011} it was shown experimentally that the flow from a swimming bacterium is well approximated by the flow generated by a force dipole. 
Above we mentioned that such a representation is not sufficient  
for our purposes to study the impact of flexible flagellum. However, for the sake of simplicity  the force dipole model is sufficient if we allow the dipole not to be straight: the line connecting points where two opposite forces are exerted is not necessarily parallel to these forces  (see Figure \ref{fig:dipoles}, right).   

\begin{figure*}[htbp]
	\begin{center}
	\begin{tabular}{cc}
		\includegraphics[width=0.475\textwidth]{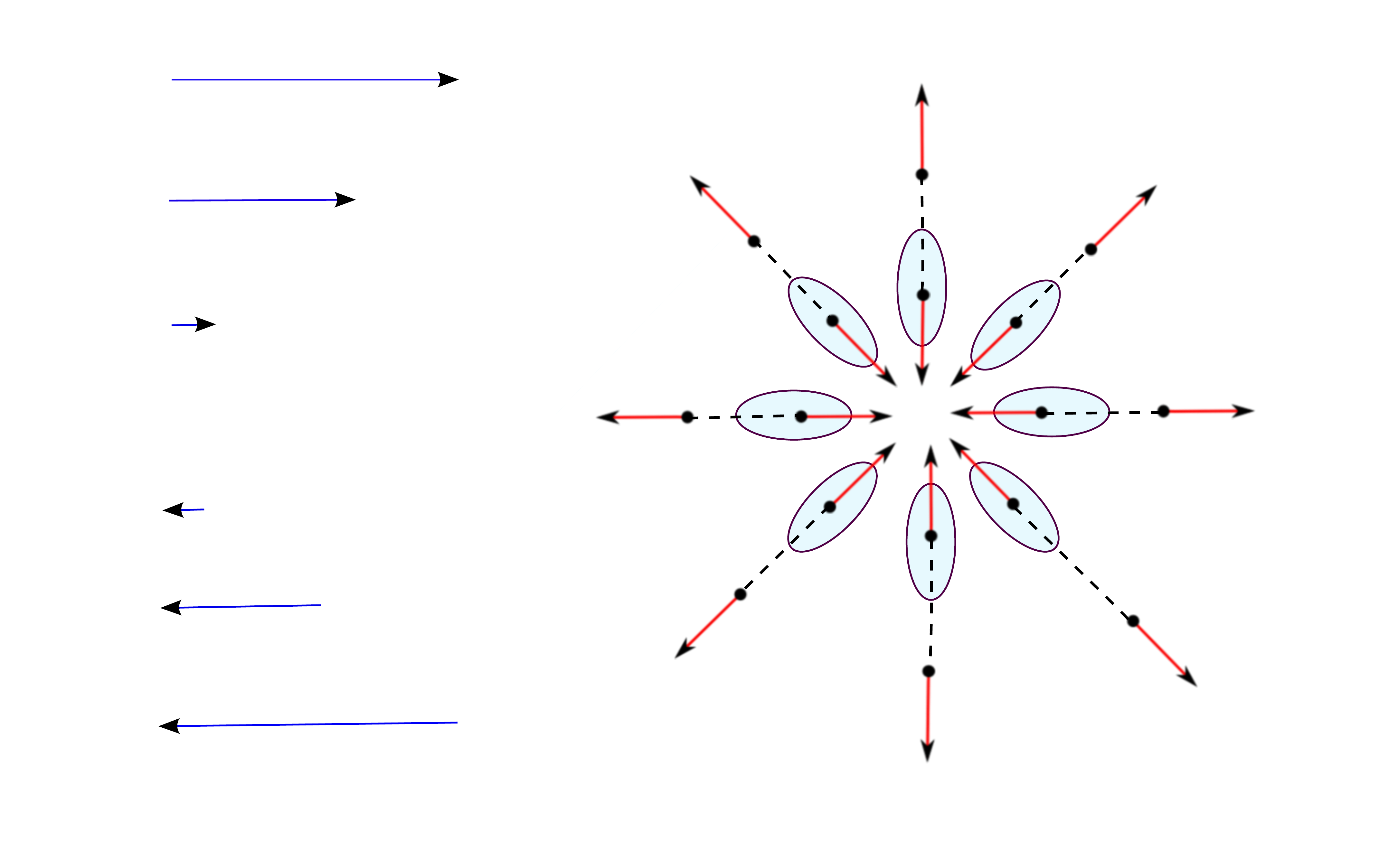}&
		\includegraphics[width=0.475\textwidth]{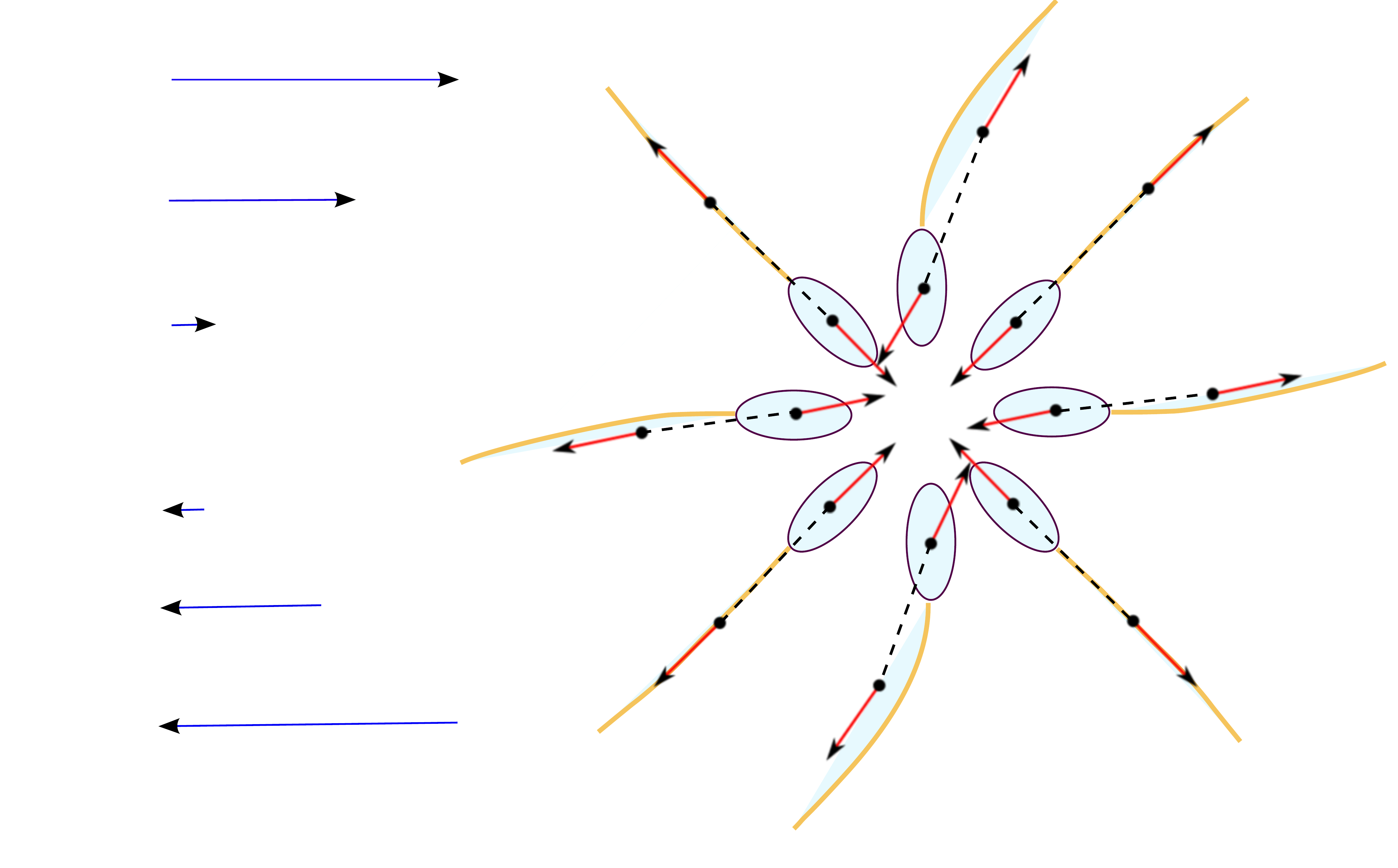} 
	\end{tabular}
	\caption{{\bf Illustration of the viscosity reduction due to flagella bending}. In red: Dipole exerted by the swimmer on the fluid, for different body orientations $i \pi/4$,  $i = 0 \dots 7$. (a) rigid flagellum. (b)  elastic flagellum.
		The shape of the flagellum is obtained from the analytic formula provided by two scale asymptotic expansion. }
	\label{fig:dipoles}
	\end{center}
\end{figure*}

The key point is that the orientation of the dipole created by each swimmer is modified by the presence of the flagellum, which, in turn, affected by the shear flow. 
 The shape of the flagellum is plotted in Figure \ref{fig:dipoles}.
Two mechanisms are responsible for the shape of the flagellum.
First,
the fluid has a direct action on the flagellum. 
This effect is prevalent when the flagellum is perpendicular to the fluid ($\theta_0 = \frac{\pi}{2}$ in Figure \ref{fig:dipoles}).
When the flagellum is oriented in the direction of the flow, this effect tends to become negligible.
Second, the rotation of the body makes one of the end of the flagellum (the end attached to the body) 
move faster than the other, which also modifies the 
shape of the flagellum.
This effect is prevalent for $\theta_0 = 0$ in Figure~\ref{fig:dipoles}.

Given that a dipole oriented a $\pi/4$ or $5\pi/4$ helps the fluid to flow, and that 
a dipole oriented a $3\pi/4$ or $7\pi/4$ prevents the fluid to flow,
the symmetry in the body orientation distribution leads to no decrease of viscosity. 
The flagellum breaks the symmetry in the dipole orientation, for example, a body oriented at $\pi/2$ (neutral for $\eta_{\text{eff}}$ without flagellum) 
creates a dipole oriented 
close to $\pi/4$ thanks to the flagellum. As a result, in average, the orientation 
helping the fluid to flow are more likely than the others.

We point out that in contrast with the work \cite{RyaHaiBerZieAra2011}, 
the propulsion contribution to $\eta_{\text{eff}}$ has no apparent singularity for the strain rate $\dot \gamma = 0$ (which is regularized by an infinitesimal rotational diffusion), and in the first approximation (rigid flagellum) it does not depend on the shear rate.
This singularity is regularized because the non dimensional parameter $\eps$ is proportional to $1/\dot \gamma$.
The bulk stress depends on how much the flagellum bends, which is in the first approximation directly proportional to the bulk rate of strain.
As a result, 
\textcolor{black}{the shear rate modifies linearly the propulsion stress, then their ratio is constant in $\dot\gamma$:  as $\dot\gamma \to 0$ (and, thus, $\varepsilon\to 0$), the propulsion contribution to the effective viscosity  equals to $-\Phi\dfrac{L^6\comMisha{\zeta_b} F_p}{K_b}Z_{\text{prop}}(\beta,r)$  (the expression  for the non-dimensional parameter $Z_{\text{prop}}(\beta,r)$  can be found in the electronic supplementary material, section 2.5)}. Thus, the small strain rate  $\dot \gamma $ limit is well-defined, and leads to a specific value of the effective viscosity even in the absence of fluctuations. This result is  in agreement with the experiment \cite{LopGacDouAurCle2015} where a well-defined value of effective viscosity was observed for very small shear rates.

\newpage 
\section*{Conclusions}
\subsection{•}
In this work we demonstrated how flexibility of bacterial flagella affects macroscopic properties of the suspension of microswimmers. 
We found that flagella bending may lead to a decrease of the effective viscosity in the absence of random reorientations. This effect is amplified  with the increase in the viscosity 
of suspending fluid since many bacteria often increase their  propulsion force  \cite{GreCan1977,SchDoe1974,Kel1974}.
Moreover, we show that flagella buckling may assist bacteria to escape entrapment at the wall. Our findings highlight the wealth of new intriguing phenomena stemming from the flexibility of the swimmer's body that include reduction of the viscosity, escape from the wall entrapment, migration towards flow centerline and many others. 

In the course of our work we approximated helical flagella by an elastic beam with the propulsion force distributed uniformly along the beam.  Obviously, this approximation neglects intrinsic 
chirality of the flagella\comMisha{, which leads to its clockwise rotation and counter-clockwise rotation of the head}. The chirality of the flagellum can be responsible for such phenomena as rheotaxis \cite{fu2012bacterial} and circular motion near the wall \cite{diluzio2005escherichia}.  Incorporating flagella chirality into our analysis would be desirable, but technically challenging. We anticipate that the torques arising from the helical shape of the flagellum are negligible compared to the bending stresses considered here, and, thus, do not affect the phenomena considered in this work (see the electronic supplementary material, section 3).


\begin{methods} 
\subsection{•}

  We use an extension of the model of the flagellated swimmer from Ref.\cite{TouKirBerAra2015}. 
  The model is two-dimensional and describes the swimmer as a rigid body of an ellipsoidal shape with an attached elastic beam representing the flagellum (see Supplementary Figure 1). The orientation of the swimmer $\theta_0$ is defined as the orientation of the principal axis of the body with respect to horizontal, and the equation for $\theta_0$ is derived from the torque balance:
  \begin{equation}
  \label{body_torque_balance}
   \zeta_r \omega = \text{T}_{\text{shear}} + \text{T}_{\text{flagellum}} + 
   \text{T}_{\text{external}}.
    \end{equation}    
 The equation reads as follows. The viscous torque which is linearly proportional to the angular velocity of the body $\omega=\dfrac{d\theta_0}{dt}$ is balanced by the torques coming from the shear flow, flagellum, and possible external torque (for example, due to the presence of a wall). Equation \eqref{body_torque_balance} is a modification of the well-known Jeffery equation for rotating ellipsoids in a background shear flow \cite{Jef1922,KimKar1991}. 
   
 
 The flagellum is a segment of a curve of the constant length $L$ and with arc-length parameter $s$, $0\leq s \leq L$. Unit vectors $\boldsymbol{\tau}(s)$, $\text{n}(s)$ and $\text{\bf b}(s)$ represent tangent, normal and bi-normal vectors of the curve, respectively. Within the flagellum, the viscous (drag) force is balanced by the propulsion and elastic forces:
 \begin{equation}
 \label{flagellum_force_balance}
 \text{F}_{\text{drag}}(s)=\text{F}_{\text{propulsion}}(s)+\text{F}_{\text{elastic}}(s).
 \end{equation}  
 The propulsion force is assumed to have a constant magnitude and to be always exerted in the tangent direction: $\text{F}_{\text{propulsion}}=F_p\boldsymbol{\tau}$. The elastic force is given by the internal stress $Q(s)$ which is the force exerted by the segment $[s,L]$ of the flagellum on the segment $[0,s]$. Thus, $\text{F}_{\text{elastic}}=-\dfrac{\partial Q}{\partial s}$. Elasticity of flagellum is constituted through the relation for internal torque $\text{\bf M}(s)~=~\boldsymbol{\tau}(s)~\times~ \text{F}_{\text{elastic}}(s)$:
 \begin{equation}
 \label{constitutive_relation}
 \text{\bf M}(s) = K_b \kappa(s) \text{\bf b}(s),
 \end{equation} 
  which reads that $ \text{\bf M}$ is proportional to the bending stiffness $K_b$ and to the local curvature of the flagellum $\kappa$. One end, $s=0$, is rigidly attached to the body (clamped), and another end, $s=L$, is free and the flagellum is straight there: $Q(L)=\kappa(L)=0$. The propulsion is transmitted to the body, thus, pushing the swimmer forward, through the point of junction and is present in the force balance equation for the body
  \begin{equation}
  \label{body_force_balance}
  \zeta_bV = \text{F}_{\text{propulsion}}(0)+\text{F}_{\text{elastic}}(0),
  \end{equation}  
 where $\zeta_b V$ is the drag force for the body and, by the Stokes law, it is proportional to the body velocity $V$ with a drag coefficient $\zeta_h$. 
 
 After proper rewriting, all the above equations \eqref{body_torque_balance}, \eqref{flagellum_force_balance}, \eqref{constitutive_relation}, \eqref{body_force_balance} result into a coupled system of an ordinary differential equation for the body, nonlinear elliptic partial differential equation of the second order for the tangential component of $Q$, highly nonlinear parabolic partial differential equation of the fourth order for the shape of the flagellum, and certain boundary conditions. The system is written in the electronic supplementary material, section 1.     
 
The system provides a deterministic description of an isolated swimmer, given its initial position, shape and orientation. 
Swimming can be defined \cite{AloDeSLef2009} as the ability to advance in a fluid in absence of external propulsive forces, by performing cyclic shape changes.
Specifically, we can define a low-Reynolds swimmer as an object which can modify, via an intrinsic mechanism (cyclic change of shape, rotation of an helix) the fluid velocity around itself, which leads to the creation of a propulsion force 
responsible for a net displacement of the center of mass.  For instance, {\it B. subtilis} actuate passive helical filaments (the flagella) using rotary motors embedded in the cell walls, and whose rotation gives rise to propulsion. This is made possible because the flagellum is an helix undergoing a drift across streamlines due to its chirality \cite{Mar2011}, as opposed to rods undergoing classical Jeffery orbits \cite{KimKar1991}. The non-symmetry in the chirality force breaks the scallop theorem\cite{BerAnd1973}.

In the derivation of the system, the following major simplifications were made. First, although the three-dimensionality of the helix is crucial for motion, the model we consider is two-dimensional, and we represent the 3D chiral force by a propulsion force density uniformly applied along the flagellum.
Second, the drag force acting on the swimmer is given by the fluid
velocity relative to the swimmer velocity.
This means that the local effects on the fluid in the neighborhood of the swimmer are not described. Such an approach is justified because we only consider dilute suspensions (the other swimmers are far). 
We also note here that somewhat similar model was considered in Ref. \cite{RopDreBauFerBibSto2006} to study magneto-elastic filaments.
We conclude this section with the list of parameters used in the paper, see Table \ref{table_2}.\\
\begin{table*}
\begin{center} 
	\begin{tabular}{| c | c | c |  }
		\hline
		parameter & typical value & description\\\hline
		$L$& $1.2\cdot 10^{-5}\;\text{m}$& flagellum length \\
		$\ell$&$0.5\cdot 10^{-5}\;\text{m}$ & body length (major axis of ellipse)\\
		$d$& $7\cdot 10^{-7}\;\text{m}$ & body thickness (minor axis of ellipse) \\
		$\beta$ &0.0162 & body shape parameter, $\dfrac{d^2}{\ell^2+d^2}$\\
		$\dot\gamma$ &	0.1 $\;\text{s}^{-1}$& shear rate \\
		$\eta_0$ & $10^{-3}\;\text{Pa}\;\text{s}$ & viscosity of the surrounding fluid \\
		$F_p$ & $10^{-7}\;\text{N} \;\text{m}^2$ & propulsion  force \\
		$K_b$	& $3\cdot 10^{-23}\;\text{N}\;\text{m}^2$	& flagellum bending stiffness \\
		\comMisha{$\zeta_b$} & $10^{-3}\; \text{N}\;\text{s}\;\text{m}^{-2}$ & drag coefficient \comMisha{per unit length} for the flagellum  \\ 
		$\zeta_h$ & $1.6 \cdot 10^{-8} \;\text{N}\;\text{s}\;\text{m}^{-1}$ & drag coefficient for the body  \\ 
		$\zeta_r$ & $6.7 \cdot 10^{-20} \;\text{N}\;\text{s}\;\text{m}$ & rotational drag coefficient for the body  \\ 
		$\alpha$ & 2 & drag anisotropy factor (ratio tangential/normal force \\ && needed to drag flagellum point) \\
		$k_r$ & 0.65 & $L \comMisha{\zeta_b}/\zeta_h$ (auxiliary parameter)\\
		$\comMisha{r}$ & 0.41 & $\ell/L$ \\
		$\comMisha{\eps}$ & 0.07 & $L^4 \dot \gamma \comMisha{\zeta_b}/K_b$ \\
		\bottomrule
\end{tabular}
	\caption{Main model parameters}
	\label{table_2}
	\end{center} 
	\end{table*}
\end{methods}

\begin{addendum}
\item
The work of M.T., M.P. and L.B. was supported by the NIH
 grant 1R01GM104978-01. The research of I.S.A. was supported by the US Department of Energy (DOE), Office of Science, Basic Energy Sciences (BES),
Materials Science and Engineering Division.

\item[Contributions]
M.P., M.T., L.B.  and I.S.A. designed and performed the study and  wrote the paper.
M.P. and M.T. contributed equally to this work.

\item[Competing Interests] The authors declare that they have no
competing financial interests.
 \item[Correspondence] Correspondence and requests for materials
should be addressed  Igor S Aranson ~ (email: aronson@anl.gov).
\end{addendum}
\appendix\onecolumn
 
\section { MMFS - Mathematical Model of Flagellated Swimmer }
			\begin{wrapfigure}[17]{r}{0.38\textwidth}
				\fbox{
					\begin{minipage}{0.38\textwidth}
						\begin{center}
							\includegraphics[width =0.8 \textwidth]{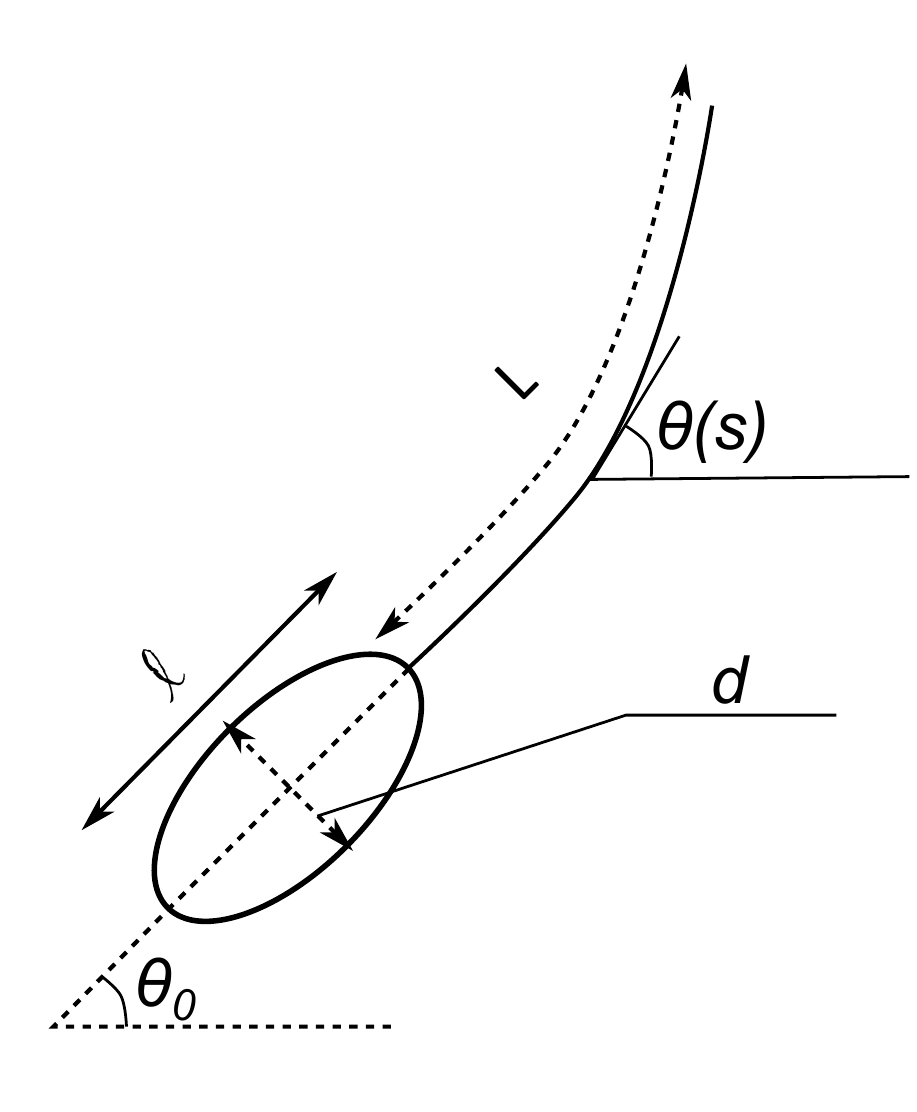}
							\caption{\small \newline Scheme of flagellated microswimmer}
							\label{fig:scheme}
						\end{center}
					\end{minipage}}
				\end{wrapfigure} 
		
		In this supplementary section, we explain how the trajectory of the swimmer and the dynamics of its flagellum are computed in the framework of MMFS. Recall that the swimmer consists of a rigid body and a flexible flagellum (see Fig.~\ref{fig:scheme}).

		The body is an ellipse centered at $X^b(t)=(x^b(t),y^b(t))$ with the major and minor axes $\ell$ and $d$, respectively. The swimmer swims in the background flow $u^{\text{BG}}$ (either the shear flow $u^{\text{BG}}(x,y)=(-\dot{\gamma}y, 0)$ or zero flow $u^{\text{BG}}=(0,0)$).  The body velocity ${V}^{b}(t)=\dfrac{\text{d}X^b(t)}{\text{d}t}$ is given by 
	\begin{equation}\label{body_move}
		V^b(t)=u^{\text{BG}}(X^b(t))+\dfrac{1}{\zeta_b}\left\{\Lambda(0)\boldsymbol{\tau}+\dfrac{1}{\alpha}N(0)\text{n}\right\}.
	\end{equation}
	Equation \eqref{body_move} reads as follows: relative velocity of the body with respect to background flow is determined by the force exerted by the flagellum (Stokes law); the parameter $\alpha$ takes into account that the ability of the flagellum to drag the body (or, the ability to affect the body velocity) in normal and tangent directions are different, $\alpha \neq  1$. Recall that $Q(s)=\Lambda(s)\boldsymbol{\tau}+N(s)\text{n}$ is the elastic stress, that is the force exerted by the segment $[s,L]$ of the flagellum on the segment $[0,s]$.

	The  flagellum location $X(s,t)=(x(s,t),y(s,t))$ 
	as a function of the body location
	$X^b(t) = (x^b(t),y^b(t))$ and flagellum orientation $\theta(s,t)$ is given by the geometrical relations
	\begin{equation}
		\left\{
		\begin{aligned}
			x(s,t)& = x(0,t)+  \int_0^s  \cos (\theta(z,t))dz,  \quad
			y(s,t) = y(0,t) + \int_0^s  \sin (\theta(z,t)) dz,  \\
			x(0,t) &= x^b(t)+\dfrac{\ell}{2} \cos (\theta_0(t)), \quad
			y(0,t)  = y^b(t)+\dfrac{\ell}{2} \sin (\theta_0(t)) .
		\end{aligned}
		\right.
	\end{equation}

	In order to find $\theta(s,t)$, $\theta_0(t)$, $\Lambda(s,t)$ and $N(s,t)$, the following PDE/ODE system is considered. It consists of an ODE for body orientation angle $\theta_0(t)$, a parabolic forth-order PDE for flagellum orientation angle $\theta(s,t)$, and an elliptic second order PDE for the tangential elastic stress $\Lambda(s,t)$. Independent variables are $t>0$ and $0< s < L$.
	\begin{equation}
		\label{theta0}
		\frac{d\theta_{0}}{dt}=-\dot\gamma \left((1-\beta)\sin^{2}\theta_{0}+\beta\cos^{2}\theta_{0}\right)+\frac{\ell }{2\zeta_r}N(0,t),
	\end{equation}

	\begin{equation}
		\label{theta}
		{\dfrac{\partial\theta}{\partial{t}}=
			-\dfrac{K_b}{\a \zeta_f}\frac{\partial^{4}\theta}{\partial{s}^{4}} +\dfrac{1}{\zeta_f}\left(\dfrac{1}{\a}{\Lambda}+ K_{b}\left(\frac{\partial\theta}{\partial{s}}\right)^{2}\right)\frac{\partial^{2}\theta}{\partial{s}^{2}} 
			+\dfrac{1}{\zeta_f}\left(\dfrac{\a+1}{\a}\frac{\partial{\Lambda}}{\partial{s}}+  F_p\right)\frac{\partial\theta}{\partial{s}}
			-\dot\gamma\sin^{2}\left(\theta\right)},
	\end{equation}
	\begin{equation}
		\begin{aligned}
			\label{Lambda}
			{\frac{\p^{2}\Lambda}{\p s^{2}}=\frac{1}{\a}  \left(\frac{\partial\theta}{\partial{s}}\right)^{2} \Lambda  
				-K_b\left(\frac{\partial^{2}\theta}{\partial{s}^{2}}\right)^{2}-\frac{\dot{\gamma}\zeta_f}{2}\sin(2\theta)
				-\frac{(\a+1)}{\a}K_b\frac{\p^{3}\theta}{\p s^{3}}\frac{\p \theta}{\p s}}. 
		\end{aligned}
	\end{equation}
	The system is supplemented with an expression for the normal component of internal stress 
	\begin{equation*}
		N=-K_b\dfrac{\p^2 \theta}{\p s^2},\quad 0\leq s \leq L, \; t>0.
	\end{equation*}
	\noindent The system is also endowed with boundary conditions at $s=0$ (interface body/flagellum):
	\begin{equation}\label{N0}
		\theta|_{s=0}=\theta_{0},
	\end{equation}
	\begin{equation}\label{krLambda}
		\dfrac{1}{\zeta_b}\Lambda|_{s=0} =\dfrac{\a \ell }{ 4 }\sin(2\theta_0) +\dfrac{1}{\zeta_f}\left[\frac{\partial{\Lambda}}{\partial{s}}|_{s=0}+ F_p +K_b\dfrac{\p \theta }{\p s}|_{s=0} {\dfrac{\partial^2 \theta}{\partial s^2}}|_{s=0}\right],
	\end{equation}
	\begin{equation}\label{krN0}
		-\left(\dfrac{1}{\alpha \zeta_b}+\dfrac{\ell^2}{4\zeta_r}\right) K_b{\dfrac{\partial^2\theta}{\partial s^2}}|_{s=0} =\dfrac{\beta \dot{\gamma}\ell}{2}\cos(2\theta_0)
		+\dfrac{1}{\alpha \zeta_f}\left[- K_b \frac{\p^3 \theta}{\p s^3}|_{s=0}+\dfrac{\p \theta}{\p s}\Lambda|_{s=0}\right],
	\end{equation}
	
	\medskip 
	
	\noindent and at $s=L$ 
	(free end of the flagellum): 
	\begin{equation}\label{bc_s11}
		\frac{\partial\theta}{\partial {s}}|_{s=L}=\frac{\partial^{2}\theta}{\partial {s}^{2}}|_{s=L}={\Lambda}|_{s=L}=0.
	\end{equation}
	
	\noindent{\bf Remark:} We explain now all the drag coefficients appearing in the system above:
	$\zeta_b$,  $\zeta_f$, $\zeta_r$, and $\alpha$. To drag the ellipsoidal body with the given velocity $V$  and angular velocity $\omega$ one need to exert the force $\zeta_b V$ and the torque $\zeta_r \omega$, respectively. It can be shown for ellipsoids that the drag coefficients $\zeta_b$ and $\zeta_r$ are related through the following expression: 
	\begin{equation}
		\zeta_r=\dfrac{\ell^2}{6}\zeta_b.
	\end{equation}
	To drag an infinitesimal (small) piece of the flagellum of length $\Delta s$ with velocity $V_{\boldsymbol{\tau}}\boldsymbol{\tau}+V_{\text{n}}\text{n}$ one needs to exert the drag force $\zeta_f\Delta s\left(\alpha V_{\boldsymbol{\tau}}\boldsymbol{\tau}+V_{\text{n}}\text{n}\right) $. As it is was mentioned above,  the parameter $\alpha$ takes into account that drag coefficients in normal and tangent directions are different for the flagellum.

	\section {Results of Two Scale Asymptotic Expansions}
	
	\subsection{Original PDE system non-dimensionalized}
	
	After the non-dimensionalization
	\begin{equation*}
		\begin{aligned}
			&	{\tilde{s}}=\frac{s}{L},\quad{\tilde{t}}=\gamma t,\quad {\tilde{\Lambda}}=\frac{\Lambda}{\zeta_{f}\gamma L^{2}}
		\end{aligned}
	\end{equation*}     
	we obtain the following PDE system of MMFS with $t>0$  and $0<s<1$:
	\begin{empheq}[left=\empheqlbrace]{align}
		\dfrac{d\theta_{0}}{d{t}}&=-\left((1-\beta)\sin^{2}\theta_{0}+\beta\cos^{2}\theta_{0}\right)+\frac{3k_r}{r} {N}_{0},\label{theta0_si}\\
		\dfrac{\partial\theta}{\partial{t}}&=
		-\dfrac{1}{\eps\a}\frac{\partial^{4}\theta}{\partial{s}^{4}} +\left(\dfrac{1}{\a}{\Lambda}+ \dfrac{1}{\eps}\left(\frac{\partial\theta}{\partial{s}}\right)^{2}\right)\frac{\partial^{2}\theta}{\partial{s}^{2}} 
		+\left(\dfrac{\a+1}{\a}\frac{\partial{\Lambda}}{\partial{s}}+  f_p\right)\frac{\partial\theta}{\partial{s}}
		-\sin^{2}\left(\theta\right),\label{theta_si}\\
		\frac{\p^{2}\Lambda}{\p s^{2}}&=\frac{1}{\a}  \left(\frac{\partial\theta}{\partial{s}}\right)^{2} \Lambda  
		-\dfrac{1}{\eps}\left(\frac{\partial^{2}\theta}{\partial{s}^{2}}\right)^{2}-\frac{1}{2}\sin(2\theta)
		-\frac{(\a+1)}{\a \eps}\frac{\p^{3}\theta}{\p s^{3}}\frac{\p \theta}{\p s}, \label{lambda_si}
	\end{empheq}
	where we dropped tildes in notations for $s$, $t$, and $\Lambda$, as well as introduce the additional parameters:  
	\begin{equation*}
		\begin{aligned}
			\varepsilon= \dfrac{\zeta_f \gamma L^4}{K_b},\qquad f_p= \dfrac{F_p}{\zeta_f \gamma L},\qquad k_r=\frac{L\zeta_{f}}{\zeta_{b}}, \qquad r=\dfrac{\ell}{L}.
		\end{aligned}
	\end{equation*}
	In what follows, we obtain asymptotic formulas in the limit $\varepsilon \to 0$. 
	
	\medskip 
	
	\noindent Equations at $s=0$:

	\begin{empheq}[left=\empheqlbrace]{align}
		N_0&= -\dfrac{1}{\eps}\;\frac{\partial^{2}\theta}{\partial{s}^{2}}|_{s=0}, \quad 	\theta|_{s=0}=\theta_{0}, \label{theta_bcs0}\\
		k_r{\Lambda}|_{s=0} &=\dfrac{\a \,r }{ 4 }\sin(2\theta_0) +\frac{\partial{\Lambda}}{\partial{s}}|_{s=0}+ f_p - \dfrac{\p \theta }{\p s}|_{s=0}N_0,	\label{lambda_bcs0}\\
		\sigma k_r{N}_{0} &=\beta\dfrac{\a\, r}{2}\cos(2\theta_0)
		- \dfrac{1}{\eps}\; \frac{\p^3 \theta}{\p s^3}|_{s=0}+\dfrac{\p \theta}{\p s}|_{s=0}\Lambda|_{s=0}.\label{n_bcs0}		
	\end{empheq}
	Here $\sigma=1+\dfrac{3\alpha}{2}$.

	\noindent Equations at $s=1$: 
	\begin{equation}\label{bc_s1}
		\frac{\partial\theta}{\partial {s}}|_{s=1}=\frac{\partial^{2}\theta}{\partial {s}^{2}}|_{s=1}={\Lambda}|_{s=1}=0.
	\end{equation}

	\subsection{Multiscale expansion}
	
	\begin{equation*}
		\left\{
		\begin{aligned}
			\theta(s,t,\tau)& = \theta^0(s,t,\tau)+ \eps \theta^1(s,t,\tau)+...,   \\ \theta_0(t,\tau)& = \theta_0^0(t,\tau) + \eps \theta_0^1(t,\tau) +..,\\
			\Lambda(s,t,\tau)  &= \Lambda^0(s,t,\tau)  + \eps \Lambda^1(s,t,\tau) +..., \\ 
			N_0(t,\tau) &= N_0^0(t,\tau) + \eps N_0^1(t,\tau)+...,
		\end{aligned}
		\right.
	\end{equation*}
	Here $\tau= \eps t$ (slow time). 
	

	\subsection{Modified Jeffery Equation for $\theta_0^0$}
	\noindent
	The equation for $\theta_0^0$ is obtained by collecting all terms at level $\eps^0$ in \eqref{theta0_si}: 
	\begin{equation}\label{mod_jeffery_unknown}
		\dfrac{\p\theta_{0}^0}{\p{t}}=-\left((1-\beta)\sin^{2}\theta_{0}^0+\beta\cos^{2}\theta_{0}^0\right)+\frac{3k_r}{r} {N}_{0}^0
	\end{equation}
	This is not a closed equation for $\theta_0^0$ because of the unknown term $\dfrac{3k_r}{r}N_0^0$ in the equation \eqref{mod_jeffery_unknown}. In order to find this term (in terms of $\theta_0^0$) first note that due to the first equation in \eqref{theta_bcs0} we have 
	\begin{equation}
		\label{what_is_N_0}
		N_0^0=-\dfrac{\p^2\theta^1}{\p s^2}|_{s=0}.
	\end{equation} 
	
	Next, expanding equation \eqref{theta_si} and collecting all terms at $\eps^{-1}$ we get 
	\begin{equation}
		0=-\frac{1}{\alpha} \frac{\partial^4\theta^0}{\partial s^4} + \left(\frac{\partial \theta^0}{\partial s}\right)^2 \frac{\partial^2 \theta^0}{\partial s^2}, \quad 0<s<1,\\
	\end{equation}
	To write boundary conditions $\theta^0$ at $s=0$ collect terms at $\eps^{-1}$ in the first equation in \eqref{theta_bcs0} and \eqref{n_bcs0}:
	\begin{equation}
		\frac{\partial^2 \theta^0}{\partial s^2}|_{s=0}=\frac{\partial^3 \theta^0}{\partial s^3}|_{s=0}=0.
	\end{equation}
	Equations in \eqref{bc_s1} give boundary conditions for $\theta^0$ at $s=1$:
	\begin{equation}
		\frac{\partial \theta^0}{\partial s}|_{s=1}=\frac{\partial^2 \theta^0}{\partial s^2}|_{s=1}=0.
	\end{equation}
	Thus, $\theta^0$ does not depend on $s$ and $\theta^0(s,t,\tau)=\theta^0_0(t,\tau)$:
	\begin{equation}
		\frac{\partial^i \theta^0}{\partial s^i}\equiv 0, \;\;i=1,2,3,... \label{theta0_is_0}
	\end{equation}
	
	\noindent 
	Due to \eqref{theta0_is_0}, the equation for $\theta^1$ which is obtained by collecting all terms at $\eps^0$ in can be written as 
	\begin{equation}\label{p4theta1}
		\frac{\partial^4 \theta^1}{\partial s^4}=\mathop{\underbrace{-\alpha \left(\frac{\partial\theta^0_0}{\partial t}+\sin^2 \theta_0^0\right)}}_{=:\chi}
	\end{equation}
	Thus, taking into account \eqref{bc_s1} we obtain 
	\begin{eqnarray}
		\frac{\partial^3 \theta^1}{\partial s^3}&=& \chi\cdot(s-1)+C_1,\label{p3theta1}\\
		\frac{\partial^2 \theta^1}{\partial s^2}&=& \frac{1}{2}\chi\cdot(s-1)^2+C_1\cdot(s-1),\label{p2theta1}
	\end{eqnarray}
	Here $C_1$ may depend on $t$ and $\tau$, but not on $s$. In order to find $C_1$, use \eqref{what_is_N_0}, \eqref{theta0_is_0} to substitute \eqref{p3theta1} and \eqref{p2theta1} with $s=0$ into \eqref{n_bcs0} at level $\eps^0$:
	\begin{equation}\label{c1_equals_to}
		C_1=\frac{\sigma k_r+2}{2(\sigma k_r +1)}\chi+ \frac{\beta\alpha \,r }{2(\sigma k_r+1)}\cos(2\theta_0^0).
	\end{equation}
	Now we are in position to find the unknown term in the equation \eqref{mod_jeffery_unknown}: 
	\begin{eqnarray}
		\frac{3k_r}{r} N_0^0&=&-\frac{3k_r}{r} \frac{\partial^2 \theta^1}{\partial s^2}|_{s=0}=\frac{3k_r}{2r}(-\chi+2C_1)\nonumber\\ 
		&=& \frac{3k_r}{2r(\sigma k_r+1)} (\chi+\beta \alpha\, r\,\cos (2\theta_0^0))\nonumber\\
		&=& -\frac{3\alpha k_r}{2r(\sigma k_r+1)}\left\{\frac{\partial \theta_0^0}{\partial t}+\sin^2{\theta_0^0}\right\}+\frac{3\beta\alpha k_r}{2(\sigma k_r+1)}\cos 2\theta_0^0. \label{N00_calc}
	\end{eqnarray}
	In order substitute \eqref{N00_calc} into \eqref{mod_jeffery_unknown} we note that due to a simple trigonometric identity
	\begin{equation}\label{trig_identity}
		-\sin^2 \theta_0^0 - B\cos(2\theta_0^0)=-(1-B)\sin^2\theta_0^0-B\cos^2 \theta_0^0,
	\end{equation}
	with $B=\beta$, the equation \eqref{mod_jeffery_unknown} can be written as follows:
	\begin{equation}\label{mod_jeffery_unknown_convenient}
		\frac{\partial \theta_0^0}{\partial t}+\sin^2{\theta_0^0}=-\beta \cos 2\theta_0^0 +\frac{3k_r}{r} N_0^0.
	\end{equation}
	Use \eqref{N00_calc} to write \eqref{mod_jeffery_unknown_convenient} in the form 
	\begin{equation*}
		\left[1+\frac{3\alpha k_r}{2 r (\sigma k_r+1)}\right]\left\{\frac{\partial \theta_0^0}{\partial t}+\sin^2 \theta_0^0\right\}=-\beta \left[1-\frac{3\alpha k_r}{2(\sigma k_r+1)}\right]\cos(2\theta_0^0),
	\end{equation*} 
	or if one divides by $1+{3\alpha k_r}/{(2r(\sigma k_r+1))}$ and uses identity \eqref{trig_identity} for 
	\begin{equation}\label{new_beta}
		B=b:=\beta r \dfrac{2\sigma k_r +2 - 3\alpha k_r }{2r\sigma k_r +2r +3 \alpha k_r}, 
	\end{equation}
	then
	\begin{equation}
		\label{eq_for_theta_00_useful}
		\frac{\partial \theta_0^0}{\partial t}+\sin^2 (\theta_0^0)=-b\cos (2\theta_0^0),
	\end{equation}
	or, equivalently, 
	\begin{equation}\label{eq_for_theta_00}
		\frac{\partial \theta_0^0}{\partial t}=-(1-b)\sin^2\theta_0^0-b\cos^2 \theta_0^0.
	\end{equation}
	The equation \eqref{eq_for_theta_00} is of the form of the Jeffery equation, and the main conclusion here is that in the limit of the rigid flagellum $K_b\to \infty$ (equivalently, $\eps \to 0$), the swimmer with the body shape parameter $\beta$ behaves as the ellipse with no flagellum and with the shape parameter $b$ defined in \eqref{new_beta} in place of $\beta$. 
	
	We note that for typical values of parameters $\sigma$, $k_r$, and $\alpha$, the parameter $b$ introduced in \eqref{new_beta} can be computed by 
	\begin{equation*}
		b=\dfrac{r\beta}{1+2r}.
	\end{equation*}
	
	\subsection{Asymptotic formula for elastic stress $Q$}
	
	In this subsection we find asymptotic formula for elastic stress $Q^0(t,\tau,s)$, whose the normal and tangential components are $N^0$ and $\Lambda^0$, respectively. The super-index $0$ means that we search for values as $\eps \to 0$. 
	
	To find $N^0$ we first note that the equality \eqref{what_is_N_0} holds for all $0\leq s \leq 1$ (not only for $s=0$ as in \eqref{what_is_N_0}):
	\begin{equation}
		\label{what_is_N_0_for_all_s}
		N^0=-\dfrac{\p^2\theta^1}{\p s^2}, \quad 0\leq s \leq 1.
	\end{equation} 
	
	From \eqref{p4theta1}, \eqref{c1_equals_to} and \eqref{eq_for_theta_00_useful} we can easily get 
	\begin{equation}
		\label{def_of_f_and_c1}
		\chi= \alpha b \cos (2\theta_0)\quad \text{and} \quad C_1 = \sigma_1 \cos (2\theta_0^0),
	\end{equation}
	where 
	\begin{equation}
		\label{def_of_sigma_1}
		\sigma_1:= \frac{\alpha b(\sigma k_r + 2) + \alpha \beta r}{2(\sigma k_r+1)} 
	\end{equation}
	Thus, from \eqref{p2theta1}, \eqref{what_is_N_0_for_all_s} and \eqref{def_of_f_and_c1} it follows that 
	\begin{equation}
		N_0= -\left(\frac{\alpha b}{2}(s-1)^2+\sigma_1(s-1)\right)\cos (2\theta_0^0).
	\end{equation}
	
	In order to find $\Lambda^0$, we collect all terms at level $\eps^0$ in the equation \eqref{lambda_si} using \eqref{theta0_is_0}:
	\begin{equation}
		\frac{\partial^2 \Lambda^0}{\partial s^2}=-\frac{1}{2}\sin(2\theta_0^0).
	\end{equation}
	In view of \eqref{bc_s1} (a boundary condition at $s=1$) and that $\theta_0^0$ is independent from $s$ we have 
	\begin{equation}
		\Lambda^0 = -\frac{1}{4}(s-1)^2 \sin 2\theta_0^0  + C_2 (s-1). 
	\end{equation}
	In order to find $C_2$ use \eqref{lambda_bcs0}:
	\begin{equation}
		-\frac{k_r}{4}\sin 2\theta_0^0 -k_rC_2= \frac{\alpha \; r}{4}\sin 2\theta_0^0 + \frac{1}{2}\sin 2\theta_0^0 +C_2 +f_p. 
	\end{equation}
	Thus, 
	\begin{equation}
		C_2 = -\sigma_2\sin 2 \theta_0^0 -\frac{F_p}{1+k_r}, \quad \text{where}\quad \sigma_2:=\frac{k_r+\alpha r+ 2 }{4(1+k_r)}.
	\end{equation}
	In particular, 
	\begin{equation}\label{formula_for_Lambda_0}
		\Lambda^0= \left(-\frac{1}{4}(s-1)^2-\sigma_2(s-1)\right)\sin 2\theta_0^0 - \frac{f_p}{1+k_r}(s-1).
	\end{equation}
	{\it Asymptotic formulas for $N^0$ and $\Lambda^0$ in the original scaling: }
	\begin{equation}\label{formula_for_N_0_L}
		N^0=-\zeta_f\gamma \left(\frac{\alpha b}{2}(s-L)^2+L\sigma_1(s-L)\right)\cos 2\theta_0^0.
	\end{equation}
	\begin{equation}\label{formula_for_Lambda_0_L}
		\Lambda^0=-\zeta_f\gamma\left(\frac{1}{4}(s-L)^2+L\sigma_2(s-L)\right)\sin 2\theta_0^0-\frac{f_p}{1+k_r}(s-L).
	\end{equation}

	\subsection{Contribution to effective viscosity}
	
	\noindent{\it The elastic contribution.} From the Kirkwood formula (see Eqn.~(2) in the paper) using integration by parts we obtain
	\begin{equation}
		\eta_\text{elastic}=\dfrac{\Phi}{2\dot\gamma\eta_0}\; <\int_0^{L} \left(Q\cdot {\bf e}_2\right)\left(\boldsymbol{\tau}\cdot {\bf e}_1\right)+\left(Q\cdot {\bf e}_1\right)\left(\boldsymbol{\tau}\cdot {\bf e}_2\right)\;ds>_{\theta_0} 
	\end{equation}
	where ${\bf e}_1=(1,0)$, ${\bf e}_2=(0,1)$, and $<\;>_{\theta_0}$ denotes the expected value with respect to probability distribution function $P(\theta_0)$ of body orientation angles (see Results section where $P(\theta_0)$ was introduced): 
	\begin{equation}\label{def_of_average}
		<g>_{\theta_0}:=\frac{q}{2\pi}\int\limits_0^{2\pi}\frac{g(\theta_0)}{1-(1-2b)\cos (2\theta_0)}d\theta_0,
	\end{equation} 
	where $q=\sqrt{1-(1-2b)^2}$. 
	
	We search for the leading term in $\eps$ which is of the order 1. Thus, assume $\boldsymbol{\tau}=(\cos \theta, \sin \theta)\approx (\cos \theta_0, \sin \theta_0)$. Then after some straightforward calculations we obtain
	\begin{equation}\label{formula_for_eta_elastic}
		\eta_\text{elastic}=\dfrac{\Phi}{\dot\gamma\eta_0}\; <\int_0^L \Lambda^0 (s;\theta_0)\sin 2\theta_0 + N^0(s;\theta_0)\cos (2\theta_0)\;ds >_{\theta_0}
	\end{equation}
	Substituting \eqref{formula_for_Lambda_0_L} and \eqref{formula_for_N_0_L} (with $\theta_0$ instead of $\theta_0^0$) into \eqref{formula_for_eta_elastic} we get 
	\begin{equation}
		\eta_{\text{elastic}}=\Phi \dfrac{L^3}{\eta_0} Z_{\text{elastic}}(\beta,r),
	\end{equation}   
	where 
	\begin{equation} 
		\begin{aligned}
			Z_{\text{elastic}}(\beta,r)&
			= \zeta_f\dfrac{r \b}{12(1+2r -2r\b)}
			\Big\{ (2r+5) \sqrt{1+2r -r \b} \Big(2r+1 -2 \sqrt{r\b}\sqrt{1+2r-\b r}\Big)\\
			&+3 \sqrt{r\b}(3r+2)(2r+1 -2 \sqrt{r\b}) \Big\}.
		\end{aligned}
	\end{equation}

	\medskip 
	
	\noindent{\it The propulsion contribution.} The Kirkwood formula for propulsive contribution has the form 
	\begin{eqnarray}
		\eta_{\text{propulsion}}&=&\dfrac{\Phi}{2\dot\gamma \eta_0} < -F_p \int_0^L\left(\boldsymbol{\tau}\cdot {\bf e}_1\right)\left(y(s)-y(0)\right) +\left(\boldsymbol{\tau}\cdot {\bf e}_2\right)\left(x(s)-x(0)\right)\;ds>_{\theta_0}\nonumber \\
		&=& - \dfrac{F_p\Phi}{2\dot\gamma \eta_0} <\int_0^L \int_0^s \cos \theta(s) \sin \theta(w) + \sin \theta(s) \cos \theta(w) \;dw \;ds >_{\theta_0} \label{def_of_ev_propulsion}
	\end{eqnarray}
	Substituting expansion $\theta=\theta^0 + \eps \theta^1 + ...$ into \eqref{def_of_ev_propulsion} after some cumbersome calculations we obtain 
	\begin{equation}\label{formula_41}
		\eta_{\text{propulsion}}=- \dfrac{F_p\Phi}{2\dot\gamma \eta_0} <\frac{L^2}{2}\sin (2\theta_0)+\eps L \cos (2\theta_0)\int_0^L\theta^1(s)\;ds>_{\theta_0} + o(\eps).
	\end{equation} 
	The first term has zero contribution after averaging with respect to $\theta_0$, and the equality \eqref{formula_41} becomes
	\begin{equation}
		\eta_{\text{propulsion}}= - \eps \dfrac{\Phi F_p L^2}{\eta_0 \dot\gamma} Z_{\text{prop}} (\beta,r),
	\end{equation}
	where
	\begin{equation} 
		\begin{aligned}
			Z_{\text{prop}}(\beta,r)&= \b \dfrac{r(17+10r)(2r+1-q\sqrt{2r+1})}{120(2r+1-\b r)^2}.
		\end{aligned}
	\end{equation}

\bibliography{flagellum}

\end{document}